\documentclass[12pt]{amsart}
\usepackage{amsfonts,color,latexsym,amssymb,amsmath,epsfig,rotating,array}

\newtheorem{dfn}{Definition}[section]
\newtheorem{rem}[dfn]{Remark} 
\newtheorem{prop}[dfn]{Proposition}

\newtheorem{thm}[dfn]{Theorem} 
 
\newtheorem{lemma}[dfn]{Lemma}
\newtheorem{cor}[dfn]{Corollary}

\newtheorem{ex}[dfn]{Example}
\definecolor{purple}{rgb}{.5,0,.5}
\definecolor{red}{rgb}{.6,0,0} 
\definecolor{green}{rgb}{0,.5,0} 

\renewcommand{\qed}{$\blacksquare$}
\newcommand{\Log}{\mathrm{Log}}

\makeatletter
\renewcommand*\env@matrix[1][c]{\hskip -\arraycolsep
  \let\@ifnextchar\new@ifnextchar
  \array{*\c@MaxMatrixCols #1}}
\makeatother

\newlength{\burg}
\newlength{\koi}
\newlength{\sma}
\newlength{\jmr}
\newcommand{\inte}{\mathrm{Int}}
\newcommand{\relint}{\mathrm{RelInt}}

\newcommand{\thth}{^{\text{\underline{th}}}}

\newcommand{\np}{{\mathbf{NP}}}

\newcommand{\bM}{\overline{M}}
\newcommand{\bC}{\overline{\cC}}

\newcommand{\Pro}{{\mathbb{P}}}

\newcommand{\Q}{\mathbb{Q}}
\newcommand{\R}{\mathbb{R}}

\newcommand{\C}{\mathbb{C}}

\newcommand{\Z}{\mathbb{Z}}
\newcommand{\bO}{\mathbf{O}}
\newcommand{\fii}{\varphi}

\newcommand{\Zn}{\Z^n}

\newcommand{\Rn}{\R^n}

\newcommand{\Cn}{\C^n}

\newcommand{\Cs}{\C^*}

\newcommand{\rank}{\mathrm{rank}} 
\newcommand{\cA}{\mathcal{A}} 
\newcommand{\cC}{\mathcal{C}} 
 
\newcommand{\cH}{H} 
 
\newcommand{\hA}{\widehat{\cA}}

\newcommand{\dia}{$\diamond$}

\newcommand{\conv}{\mathrm{Conv}}

\newcommand{\init}{\mathrm{In}}
\newcommand{\sign}{\mathrm{sign}}
\newcommand{\codim}{\mathrm{codim }}

\author{J.\ Maurice Rojas}
\email{rojas@math.tamu.edu} 
\address{TAMU 3368, College Station, TX \ 77843-3368} 
\thanks{Partially supported by NSF grant CCF-1409020} 
\author{Korben Rusek} 
\email{korben.rusek@gmail.com} 
\address{1311 H St SE, Auburn, WA \ 98002, USA } 
\title[$\cA$-Discriminants for Complex Exponents, and Isotopy]{\mbox{}\\  
\vspace{-1.1in}$\cA$-Discriminants for Complex Exponents, and Counting 
Real Isotopy Types }   

\begin{document}

\begin{abstract} 
We extend the definition of $\cA$-discriminant varieties,   
and Kapranov's parametri- 
zation of $\cA$-discriminant varieties, 
to complex exponents. As an application, we study the 
special case where $\cA$ is a fixed real $n\times (n+3)$ matrix 
whose columns form the spectrum of an $n$-variate exponential 
sum $g$ with fixed signed vector for its coefficients: 
We prove that the number of possible isotopy types for the 
real zero set of $g$ is $O(n^2)$. The best previous upper bound was 
$2^{O(n^4)}$. Along the way, we also show that the singular 
loci of our generalized $\cA$-discriminants are images of 
low-degree algebraic sets under certain analytic maps. 
\end{abstract} 

\maketitle 

\section{Introduction} 
Classifying families of real algebraic curves up to deformation
is part of Hilbert's 16$\thth$ Problem and, for curves of high
degree, remains a challenging open problem (see, e.g., \cite{kaloshin,
orevkov}). While the number of possible isotopy types is now known to
be exponential in the degree, degree is but one measure of the
complexity of a polynomial. Here, in the spirit of Khovanski's Fewnomial 
Theory \cite{few}, we will see an example of a different
measure of complexity of zero sets (and a different class of families) 
for which the number of isotopy 
types admits a polynomial upper bound. The techniques we use to 
achieve these bounds may be of independent interest. 

In particular, we take a first step toward extending the theory of
$\cA$-discriminants \cite{gkz94}, and Kapranov's parametrization
of $\cA$-discriminant varieties \cite{kapranov}, to a broader family of
analytic functions. As an application, we prove a
quadratic upper bound on the number of isotopy types of real zero sets of 
certain $n$-variate exponential sums, in a setting where the best previous 
bounds were strongly exponential in $n$. In a companion paper 
\cite{fnr}, we apply the techniques developed here to give new polynomial 
upper bounds on the number of connected components of such real zero sets. 

More precisely, let $\cA\!=\![a_{i,j}]\!\in\!\C^{n\times t}$ have 
distinct columns, with $j\thth$ column $a_j$, 
$y\!:=\!(y_1,\ldots,y_n)\!\in\!\Cn$, 
$a_j\cdot y\!:=\!\sum^n_{i=1}a_{i,j}y_i$, 
and consider exponential sums of the form \[ g(y)\!:=\!\sum^t_{j=1} 
c_j e^{a_j\cdot y},\] where $c_j\!\in\!\C\setminus\{0\}$ for all $j$. 
We call $g$ an {\em $n$-variate exponential $t$-sum},  
$\cA$ the {\em spectrum} of $g$, and we let $Z_\R(g)$ (resp.\ $Z_\C(g)$) 
denote the set of roots of $g$ in $\Rn$ (resp.\ $\Cn$). In particular, 
we call $g$ a {\em real} exponential sum when both $\cA\!\in\!\R^{n\times t}$ 
and all the coefficients $c_j$ of $g$ are real. Also, given any two 
subsets $X,Y\!\subseteq\!\Rn$, an {\em isotopy from $X$ to $Y$ (ambient 
in $\Rn$)} is a continuous map $I : [0,1]\times \Rn \longrightarrow 
\Rn$ satisfying (1) $I(t,\cdot):\Rn\longrightarrow \Rn$ is a homeomorphism 
for each $t\!\in\![0,1]$, (2) $I(0,x)\!=\!x$ for all $x\!\in\!\Rn$, and 
(3) $I(1,X)\!=\!Y$.$^1$ 
\setcounter{footnote}{1}  
\footnotetext{It is easily checked that an 
isotopy from $X$ to $Y$ implies an isotopy from $Y$ to $X$ as well. 
So isotopy is in fact an equivalence relation and it makes sense to speak of 
{\em isotopy type}.}
Although our generalized $\cA$-discriminants allow 
arbitrary $\cA\!\in\!\C^{n\times t}$, our isotopy counts will focus on 
smooth $Z_\R(g)$, in the special case where $c_j\!\in\!\R\setminus\{0\}$ for 
all $j$ and we fix both $\cA\!\in\!\R^{n\times (n+3)}$ and 
$\sign(g)\!:=\!(\sign(c_1),\ldots,\sign(c_t))$ (the {\em sign vector} of 
$g$).

A consequence of our generalized $\cA$-discriminants (defined in the next 
section, and\linebreak 
parametrized explicitly in Definition \ref{dfn:first} and Theorem 
\ref{thm:gena} below) is the following new count of isotopy types:  
\begin{thm}
\label{thm:nearckt} 
Following the notation above, assume $t\!=\!n+3$, $\{a_1,\ldots,a_{n+3}\}$ 
does not lie in any affine hyperplane, and that we fix $\cA$ and $\sign(g)$. 
Then the number of possible isotopy types for a smooth $Z_\R(g)$ is no greater 
than $\frac{1}{2}n^2+\frac{3}{2}n+6$.  
\end{thm} 

\noindent 
We prove Theorem \ref{thm:nearckt} in Section \ref{sec:topo}. 
Using the abbreviation $x^u\!:=\!x^{u_1}_1 \cdots x^{u_n}_n$ when
$(u_1,\ldots,u_n)$\linebreak 
$\in\!\Zn$, we call $f(x)\!:=\!\sum^t_{j=1}c_j x^{a_j}$
an {\em $n$-variate $t$-nomial} when $\cA\!\in\!\Z^{n\times t}$.
The change of variables $x=e^y\!:=\!(e^{y_1},\ldots,e^{y_n})$ easily shows
that, when $\cA\!\in\!\Z^{n\times t}$, studying zero sets of $t$-nomials in
the positive orthant $\Rn_+$ is the same as studying zero sets of exponential
$t$-sums in $\Rn$, up to the diffeomorphism between $\Rn_+$ and $\Rn$
defined by $x\!=\!e^y$. So in this sense, dealing with exponential sums is 
a generalization of the polynomial case. 

A result of Basu and Vorobjov \cite{basu}
implies a $2^{O(t^4)}(nt)^{O(t(t+n))}$ upper bound for the number 
of isotopy types of a smooth $Z_\R(g)$ when $g$ is an $n$-variate 
$t$-sum. Although \cite[Thm.\ 1.3]{drrs} presents an $O(n^6)$ upper bound 
for the special case $\cA\!\in\!\Z^{n\times (n+3)}$ (with some additional 
restrictions), Theorem \ref{thm:nearckt} here is sharper, more general, and 
thus gives the first polynomial upper bound on the number of isotopy types in 
the setting of $n$-variate exponential $(n+3)$-sums.  

One should observe that if one allows 
$\{a_1,\ldots,a_t\}$ to lie in an affine hyperplane, then counting isotopy 
types becomes more complicated. In essence, this is because the number of 
isotopy types for smooth $Z_\R(g)$ with underlying spectrum 
$\cA\!=\![a_1,\ldots,a_t]$ depends on the dimension of the convex hull of 
$\{a_1,\ldots,a_t\}$, as well as its cardinality. In what follows, 
for any subset $S\!\subseteq\!\Rn$, 
we let $\conv S$ denote the convex hull of\footnote{i.e., smallest convex 
set containing} $S$. 
\begin{lemma}
\label{lemma:weird} 
There exist $\cA\!\in\!\Z^{1\times t}$ (resp.\ $\cA\!\in\!\Z^{2\times t}$) 
such that the number of possible isotopy types of a smooth 
$Z_\R(g)$, with $g$ having spectrum $\cA$, is $t$ (resp.\ $2^{\Omega(t)}$). 
\end{lemma}

\noindent
We prove Lemma \ref{lemma:weird} in Section \ref{sub:weird}. Letting 
$d(\cA)\!:=\!\dim\conv\{a_1,\ldots,a_t\}$, note that\linebreak  
$1\!\leq\!d(\cA)\!\leq\!\min\{n,t-1\}$ when $t\!\geq\!2$, and $Z_\R(g)$ 
is empty when $t\!=\!1$. Note also that the cases where $t-d(\cA)$ is fixed are 
{\em not} addressed by Lemma \ref{lemma:weird}, since $d(\cA)$ is fixed 
in the statement above. In particular, the number of 
possible isotopy types for smooth $Z_\R(g)$, when $\cA$ and $\sign(g)$ are 
fixed, is known to be $2$ when $t-d(\cA)\!\leq\!2$ (see, e.g., 
\cite{snc,bihan}). Our Theorem \ref{thm:nearckt} thus addresses the case 
$t-d(\cA)\!=\!3$. Let us now introduce our main theoretical tool. 

\subsection{Generalizing, and Parametrizing, $\cA$-Discriminants for Complex 
Exponents} Our isotopy count provides a motivation for generalized 
$\cA$-discriminants, since discriminants parametrize degenerate behavior, and 
different isotopy types can be obtained by deforming a zero set through a 
degenerate state. This connection is classical and well-known. See, for 
instance, Hardt's Triviality Theorem \cite{hardt} (in the semi-algebraic 
setting) or \cite[Ch.\ 11, Sec.\ 5]{gkz94} (in the real algebraic setting): 
The connected components of the complement of the real part of a discriminant 
variety describe regions in coefficient space 
(called {\em discriminant chambers}) where the
topology of the real zero set (in a suitable compactification of 
$\Rn$) of a polynomial is constant. In what follows, 
for any $S\!\subseteq\!\C^N$ we let $\overline{S}$
denote the Euclidean closure of $S$.
Our central object of study will be the following kind of discriminant 
variety associated to families of exponential sums: 
\begin{dfn} 
Let $\cA\!\in\!\C^{n\times t}$ have $j\thth$ column 
$a_j$ and assume the columns of $\cA$ are distinct. 
We then define the {\em generalized 
$\cA$-discriminant variety} to be\\ 
\mbox{}\hfill  
$\Xi_\cA\!:=\!\overline{\left\{\left.  \left. 
[c_1:\cdots:c_t]\!\in\!\Pro^{t-1}_\C
\right\backslash \{c_1\cdots c_t\!\neq\!0\}\; \; \right| \; \; \sum^t_{j=1} 
c_je^{a_j\cdot z} \text{ has a degenerate root in } 
\Cn\right\}}$. \hfill \dia 
\end{dfn} 

\noindent 
When $\cA\!\in\!\Z^{n\times t}$ our $\Xi_\cA$ agrees with the classical 
$\cA$-discriminant variety $\nabla_\cA$ \cite{gkz94}. More to the point, for 
all but a small family of $\cA\!\in\!\Zn$, $\nabla_\cA$ is an algebraic 
hypersurface with a defining 
polynomial that is often too unwieldy for computational purposes.  
So we will also need a more efficient, alternative description for $\Xi_\cA$. 
\begin{ex} 
Let $\cA\!=\!\begin{bmatrix} 0 & 1 & 0 & 4 & 1 \\ 0 & 0 & 1 & 1 & 4 
\end{bmatrix}$. A routine {\tt maple} calculation then shows us that\linebreak 
\scalebox{.87}[1]{$\Xi_\cA\!=\!\nabla_\cA$ here, and this $\nabla_\cA$ is the 
zero set of a polynomial $\Delta_\cA\!\in\!\Z[c_1,\ldots,c_5]$ satisfying 
$\Delta_\cA(1,1,1,a,b)=$}\\   
\mbox{}\hspace{.5in}
\scalebox{.65}[.65]{\vbox{
\noindent 
$41987654504771523593992227a^8b^8
+8568922617577790827960320a^8b^7+
394594247668399678957787136a^7b^8+491069384583950065193975808a^8b^6-\
971141005960243113814917120a^7b^7+1644546811048059090366627840a^6b^8+
557969223231079901560832a^9b^4+828434941582623838008508416a^8b^5-\
8896118143687124537286066176a^7b^6+4692084142913135619868721152a^6b^7+
2845499698372999866809843712a^5b^8+557969223231079901560832a^4b^9+
33392996500536631555522560a^9b^3+384254443547034707078152192a^8b^4-\
11483443502644561069909999616a^7b^5+14323107664774924348979937280a^6b^6+
13591000063033685271054909440a^5b^7+2225676679631729339955937280a^4b^8-\
25511283567328457194995712a^3b^9+20941053496075364622925824a^9b^2+
269737322421295126029533184a^8b^3-2514558123743644571580497920a^7b^4+
21319282121430982186963566592a^6b^5+18138163316374406659527671808a^5b^6+
4594348961140867552012926976a^4b^7+38288951865122947982163968a^3b^8-\
51524645931445780035403776a^2b^9+363087263602825104457728a^{10}-\
9792009640288689535844352a^9b+178810349707236426746167296a^8b^2+
1368264254117216589547831296a^7b^3+10397247952186084766590697472a^6b^4+
6930726608820725492905672704a^5b^5+2535119422553880950892134400a^4b^6+
134703665565747736152637440a^3b^7-92973237722754317832683520a^2b^8-\
2893351631835012551147520ab^9+363087263602825104457728b^{10}+
726174527205650208915456a^9-12696707749111290371506176a^8b-\
74489621423517087836405760a^7b^2+363189381895713399018356736a^6b^3-\
482236618449489680142434304a^5b^4+191290255533750888626651136a^4b^5+
50571247933680984080252928a^3b^6-31282237054780900405936128a^2b^7-\
5798049740657613386809344ab^8+726174527205650208915456b^9+
363087263602825104457728a^8-2904698108822600835661824a^7b+
10166443380879102924816384a^6b^2-20332886761758205849632768a^5b^3+
25416108452197757312040960a^4b^4-20332886761758205849632768a^3b^5+
10166443380879102924816384a^2b^6-2904698108822600835661824ab^7+
363087263602825104457728b^8 \text{. \scalebox{1.25}[1.25]{\dia}}$}}  
\end{ex} 

The {\em generalized} $\cA$-discriminant need not be the zero set of 
{\em any} polynomial function,\linebreak 
already for $\cA\!\in\!\R^{1\times 3}\setminus \Q^{1\times 3}$:
For instance, taking $\cA\!=\!\left[0,1,\sqrt{2}\right]$, it is not hard to 
check that the intersection of the $c_1\!=\!c_3$ line with $\Xi_\cA$ in 
$\Pro^2_\C$ is exactly the infinite set\\ 
\mbox{}\hfill $\left\{\left. \left[1:\frac{-\sqrt{2}}{\sqrt{2}-1}
(\sqrt{2}-1)^{1/\sqrt{2}}e^{\sqrt{-2}\pi k}:1\right] \; \right| \; k\!\in\!\Z
\right\}$.\hfill\mbox{}\\ 
So this particular $\Xi_\cA$ can't even be semi-algebraic. 

Nevertheless, the generalized discriminant $\Xi_\cA$ admits a 
concise and explicit parametrization that will be our main theoretical 
tool. Some more notation we'll need is the following. 
\begin{dfn} 
\label{dfn:B} 
For any $\cA\!\in\!\C^{n\times t}$ let $\hA\!\in\!\C^{(n+1)\times t}$
denote the matrix with first row $[1,\ldots,1]$ and bottom $n$ 
rows forming $\cA$, and let $B\!\in\!K^{t\times (t-d(\cA)-1)}$ 
be any matrix whose columns form a basis for the right nullspace of $\hA$, 
where $K$ is the minimal field containing the entries of $\cA$. 
Let $\beta_i$ denote the $i\thth$ row of $B$. Finally, let us call $\cA$ 
{\em non-defective} if we also have that (0) the columns of $\cA$ are distinct 
and (1) $\codim_\C \Xi_\cA\!=\!1$. \dia 
\end{dfn}

\noindent
It is easily checked $d(\cA)\!=\!\left(\rank\hA\right)-1$. 
In particular, the existence of a $B\!\in\!K^{t\times (t-n-1)}$ with
columns forming a basis for the right nullspace of $\hA$ is equivalent
to $\hA$ having rank $n+1$. Note also that $\cA$ can be defective when 
$d(\cA)\!=\!n$: For instance, $\cA\!=$\scalebox{.6}[.6]{$\begin{bmatrix} 
0 & 1 & 0 & 0 & 1 & 0 \\ 
0 & 0 & 1 & 0 & 0 & 1 \\ 
0 & 0 & 0 & 1 & 1 & 1 \end{bmatrix}$} has $d(\cA)\!=\!3$ but
is defective (see, e.g., \cite[Ex.\ 2.8 \& Cor.\ 3.7]{bhpr}
or \cite{ds02,dr06}). 

Once we define a particular hyperplane arrangement it will then be easy to 
write down our parametrization of $\Xi_\cA$. 
\begin{dfn} 
\label{dfn:first} 
Let $(\cdot)^\top$ denote matrix transpose and, for any  
$u\!:=\!(u_1,\ldots,u_N)$ and\linebreak 
$v\!:=\!(v_1,\ldots,v_N)$ in $\C^N$,  
let $u\odot v\!:=\!(u_1v_1,\ldots,u_Nv_N)$. 
Then, following the notation of Definition \ref{dfn:B}, 
assume $\cA$ is non-defective, set 
$\lambda\!:=\!(\lambda_1,\ldots,\lambda_{t-d(\cA)-1})$ and 
$[\lambda]\!:=\![\lambda_1:\cdots:\lambda_{t-d(\cA)-1}]$.  
We then define the (projective) hyperplane arrangement \\ 
\mbox{}\hfill 
$\cH_\cA\!:=\!\left\{[\lambda] \; \left| 
\; \; \lambda\cdot \beta_i\!=\!0\text{ for some } i\!\in\!\{1,\ldots,t\}\right. 
\right\}\subset\Pro^{t-d(\cA)-2}_\C$, \hfill\mbox{}\\ 
and define $\psi_\cA : \left(\left. \Pro^{t-d(\cA)-2}_\C\right\backslash  
\cH_\cA\right)\times \Cn \longrightarrow \Pro^{t-1}_\C$ by 
$\psi_\cA([\lambda],y)\!:=\!
\left[(\lambda B^\top)\odot e^{-yA}\right]$. \dia  
\end{dfn} 

\noindent 
While $\psi_\cA$ certainly depends on $B$, its image is easily seen to be 
independent of $B$: Simply note that, up to transposes,  
$\left\{\lambda B^\top\right\}_{\lambda\in\Pro^{t-d(\cA)-2}_\C}$ 
is the right-null space of $\hA$. More importantly, as we'll soon see, 
$\psi_\cA$ is in fact a parametrization of $\Xi_\cA$. In what follows, 
for any fixed choice of $\sigma\!\in\!\{\pm 1\}^t$, we call 
$\Pro^{t-1}_\sigma\!:=\!\{[\lambda]\!\in\!\Pro^{t-1}_\R\; | \; 
\sign(\lambda)\!=\!\pm \sigma\}$ an {\em orthant} of $\Pro^{t-1}_\R$. 
\begin{thm} 
\label{thm:gena} 
If $\cA\!\in\!\C^{n\times t}$ is non-defective then:\\ 
1. $\Xi_\cA\!=\!\overline{\psi_\cA\left(\left(\left. 
\Pro^{t-d(\cA)-2}_\C\right\backslash H_\cA\right)\times \Cn\right)}$. 
In particular, $\Xi_\cA$ is connected, 
and is the Euclidean\linebreak 
\mbox{}\hspace{.5cm}closure of an analytic hypersurface in 
$\Pro^{t-1}_\C$. 

\noindent 
2.\ If  we also have $\cA\!\in\!\R^{n\times t}$, 
and $\sigma\!\in\!\{\pm 1\}^t$, then  
$\Log|\Xi_\cA\cap \Pro^{t-1}_\sigma|$ is an $\Rn$-bundle over a base\linebreak  
\mbox{}\hspace{.4cm}of the form $\Gamma_\sigma(\cA)\cup 
Y_\sigma\!\subset\!\R^{t-d(\cA)-1}$ where\\ 
\mbox{}\hspace{1cm}a.\ $\Gamma_\sigma(\cA)$ is the closure of the union of 
$O(t)^{4t-2d(\cA)-2}$ real analytic hypersurfaces.\\  
\mbox{}\hspace{1cm}b.\ $Y_\sigma$ is a countable, locally finite union of 
real analytic varieties of codimension $\geq\!2$.  
\end{thm} 
\noindent 
We prove Theorem \ref{thm:gena} in Section \ref{sec:back}. 
The special case $\cA\!\in\!\Z^{n\times t}$ of the first assertion recovers 
the famous {\em Horn-Kapranov Uniformization}, derived by Kapranov in 
\cite{kapranov}. 

Theorem \ref{thm:gena} is the first key idea to proving Theorem 
\ref{thm:nearckt}: When $t\!=\!n+3$, we will see that the portion of 
$\Xi_\cA$ intersecting any fixed orthant is the union of 
a finite number of real analytic hypersurfaces that is polynomial in $n$ (and 
possibly some pieces of higher codimension). 
The next key idea (described in Section \ref{sec:chamber} below) 
is to reduce the dimension so that we can reduce 
to counting the number of connected components of the complement of 
a union of arcs. The final key idea is then showing that each 
such arc is convex and has a highly restricted Gauss map, thus enabling a 
further reduction to counting cells in an arrangement of few lines. This is 
described toward the end of Section \ref{sec:chamber} and in 
Section \ref{sec:topo}.  

\section{Defining Generalized $\cA$-Discriminant Contours and Chambers} 
\label{sec:chamber} 
Let us first observe some basic cases where $\Xi_\cA$, or at least its real 
part, can be described very easily.
\begin{lemma}
\label{lemma:basic} Assume $\cA\!\in\!\C^{n\times t}$ has distinct columns. 
Then:\\ 
1. If $t-d(\cA)\!=\!1$ then $\Xi_\cA\!=\!\emptyset$.\\  
2. If $\cA\!\in\!\R^{n\times t}$, $t-d(\cA)\!=\!2$, and 
$b\!=\!(\beta_1,\ldots,\beta_t)\!\in\!(\R\setminus\{0\})^t$ is a generator for 
the right\linebreak 
\mbox{}\hspace{.5cm}nullspace of $\hA$, then $[c_1:\cdots:c_t]
\!\in\!\Xi_\cA\cap\Pro^{t-1}_\R$ is equivalent to the following
condition:\linebreak 
\mbox{}\hfill 
$\prod^t_{i=1}\left|\frac{c_j}{\beta_j}\right|^{\beta_j}\!=\!1$ and
$\sign(c_1,\ldots,c_t)\!=\!\pm \sign(b)$.\hfill \mbox{}\\  
\mbox{}\hspace{.5cm}In particular, $\Xi_\cA\cap\Pro^{t-1}_\sigma$ is non-empty 
exactly when $\sign(b)\!=\!\pm \sigma$.  
\end{lemma}

\noindent
{\bf Proof:} The first assertion follows immediately after two
observations: (a) We may assume $a_1\!=\!\bO$ since $\Xi_\cA$
is invariant under translation of the point set $\{a_1,\ldots,a_{n+1}\}$
and (b) applying the invertible change of variables $y\!=\![a_2,\ldots,
a_{n+1}]^{-1}z$ reduces us to the case where $g(z)\!=\!c_1+
c_2e^{z_1}+\cdots+c_{n+1}e^{z_n}$, which clearly yields 
$\Xi_\cA\!=\!\emptyset$.

The second assertion is contained in \cite[Thm.\ 2.9]{snc}.
(The special case of integral $\cA$ was observed earlier in
\cite[Prop.\ 1.2, pg.\ 217]{gkz94} and \cite[Prop.\ 1.8, Pg.\ 274]{gkz94}.)
\qed 

\medskip 
Let us also observe a simple property of non-defective $\cA$. 
First, let us call $\cA$ {\em pyramidal} if and only if $\cA$ has a
column $a_j$ such that $\{a_1,\ldots,a_{j-1},a_{j+1},\ldots,a_{n+k}\}$
lies in a $(d(\cA)-1)$-dimensional affine subspace. 
\begin{prop}
\label{prop:pyr}
Following the preceding notation, $\cA$ is pyramidal 
if and only if $B$ has a zero row. In particular,
$\cA$ non-defective implies that $\cA$ is not pyramidal. 
\end{prop}  

\noindent 
{\bf Proof:} The first assertion follows easily upon observing 
that a column $b_i\!=\![b_{1,i},\ldots,b_{t,i}]^\top$ of $B$ 
has $b_{i,j}\!\neq\!0$ if and only if $a_j$ lies in the same 
affine hyperplane as $\{a_1,\ldots,a_{j-1},a_{j+1},\ldots,a_t\}$, and 
the affine relations defined by $b_1,\ldots,b_{t-d(\cA)-1}$ are 
linearly independent. 

To prove the final assertion, assume $\cA$ is pyramidal. 
Observing that $\Xi_\cA$ is invariant under invertible linear maps applied to 
$\cA$, and permuting coordinates if necessary, we may assume $j\!=\!t$, 
$a_t\!=\!e_n$ (the $n\thth$ standard basis vector), and 
$a_1,\ldots,a_t\!\in\!\R^{n-1}\times \{0\}$. 
In particular, we then see that $\Xi_\cA$ is simply $\Xi_{\{a_1,
\ldots,a_{t-1}\}}\times \{0\}$, and thus $\dim \Xi_\cA\!<\!t-1$. \qed 

\medskip 
Now let $\Log|\cdot| : \C^N \longrightarrow \R^N$ denote the coordinate-wise 
log-absolute value map. Since singularities in $Z_\C(g)$ are preserved under 
translation of $\cA$ and scaling the coefficients of $g$, the discriminant 
variety $\Xi_\cA$ possesses certain homogeneities that we can quotient out to 
simplify our study of discriminant chambers. Taking $\Log|\cdot|$ helps clarify 
these homogeneities and, as we'll see at the end of this 
section, also helps clarify the geometry of the underlying quotient.  
Toward this end, observe that $[1,\ldots,1]B\!=\!\bO$. 
Clearly then, for any $\ell\!\in\!(\R\setminus\{0\})^t$, 
we have that $(\Log|\eta\ell|)B$ is independent of 
$\eta$ for any nonzero $\eta\!\in\!\C\setminus\{0\}$, so we can then define 
$(\Log|[\ell]|)B\!:=\!(\Log|\ell|)B$. 
\begin{dfn}
\label{dfn:contour} 
Following the notation of Definition \ref{dfn:B}, define\\ 
\mbox{}\hfill $\xi_{\cA,B} : 
\left(\left. \Pro^{t-d(\cA)-2}_\C\right\backslash \!H_\cA\right) \
\longrightarrow \R^{t-d(\cA)-1}$ \ \ \ 
by \ \ \  $\xi_{\cA,B}([\lambda])\!:=\!\left(\Log\left|\lambda B^\top
\right|\right)B$.\hfill\mbox{}\\
(So $\xi_{\cA,B}$ is defined by multiplying a row vector by a matrix.)
We then define the {\em reduced discriminant contour}, $\Gamma(\cA,B)$, to 
be $\emptyset$ or $\overline{\xi_{\cA,B}\!\left(\left. 
\Pro^{t-d(\cA)-2}_\R\right\backslash\!H_\cA\right)}$,  
according as $\cA$ is defective or not. \dia
\end{dfn}

\noindent 
\vbox{
\noindent 
\begin{minipage}[t]{.6 \textwidth}
\vspace{0pt}
\begin{ex} \label{ex:penta} 
When $\cA:=$\scalebox{1}[.6]{$\begin{bmatrix}0 & 1 & 0 & 4 & 1\\
0 & 0 & 1 & 1 & 4\end{bmatrix}$} we are in essence considering the 
family of exponential sums $g(y)\!:=\!f\!\left(e^{y_1},e^{y_2}\right)$ 
where $f(x)\!=\!c_1+c_2x_1+c_3x_2+c_4x^4_1 x_2+c_5 x_1 x^4_2$. A suitable 
$B$ (among many others) with columns defining a basis for the right nullspace 
of $\hA$ is then $B\!\approx\; $\scalebox{.6}[.6]
{$\begin{bmatrix}[r] 0.5079 & -0.8069 & 0.1721 & 0.2267 & -0.0997\\
0.5420 & 0.1199 & -0.7974 & -0.0851 & 0.2206 \end{bmatrix}^\top$}, and 
the corresponding reduced contour $\Gamma(\cA,B)$, intersected with $[-4,4]$,  
is drawn to the right. \dia 
\end{ex} 
Note that in our preceding example, $\Xi_\cA$ is a hypersurface in 
$\Pro^4_\C$, and $Z_\C(g)$ has a singular point $y\!\in\!\C^2$ if and only if 
$Z_\C(h)$ has singular point $y+(\delta_1,\delta_2)$, where\linebreak  
\end{minipage} \hspace{.3cm} 
\begin{minipage}[t]{.2 \textwidth}  
\raisebox{-5.5cm}{\epsfig{file=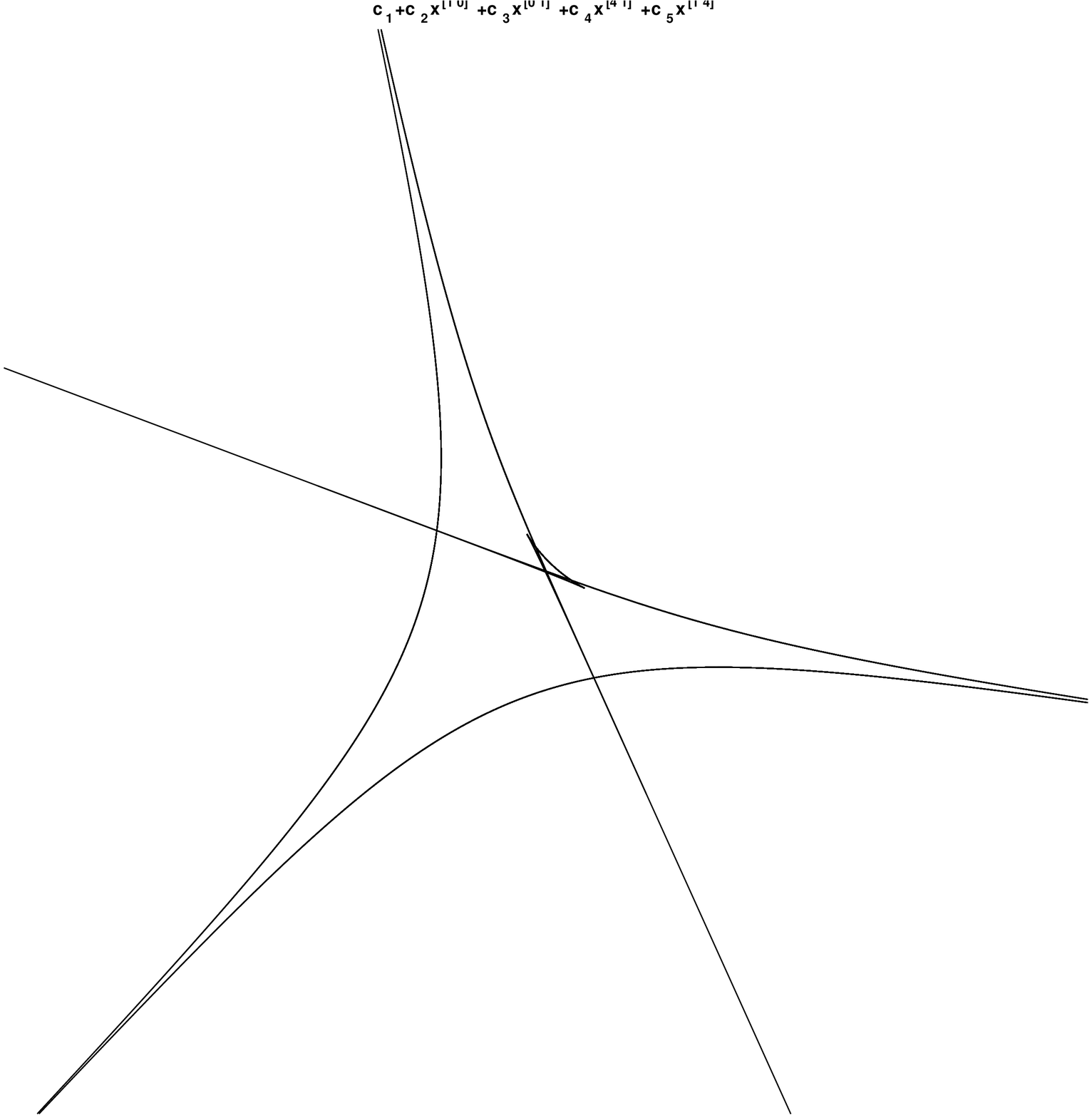,height=2.3in,clip=}}  
\end{minipage}}   

\vspace{-.3cm} 
\noindent 
$h(y)\!:=\!\alpha g(y_1-\delta_1,y_2-\delta_2)$ for some $\alpha\!\in\!\Cs$ 
and $\delta_i\!\in\!\R$.  
It is then easily checked that $\Log|\Xi_\cA\cap\Pro^4_\R|$ is a $2$-plane 
bundle over the curve in $\R^2$ drawn above. 

An obvious issue behind taking $\Log|\cdot|$ of the real part of 
$\Xi_\cA$ is that we lose information about coefficient signs. 
So we refine the notion of discriminant contour as follows: 
\begin{dfn}
\label{dfn:chamber}
Suppose $\cA\!\in\!\R^{n\times t}$ is non-defective and   
$\sigma\!=\!(\sigma_1,\ldots,\sigma_{t})\!\in\!\{
\pm 1\}^{t}$. We then define the {\em signed} reduced discriminant contour, 
$\Gamma_\sigma(\cA,B)$, to be $\emptyset$ or \\ 
\mbox{}\hfill $\overline{\left\{\left.\xi_{\cA,B}([\lambda]) \; \right| \;
\mathrm{sign}\left(\lambda B^\top\right)\!=\!\pm \sigma  \ , \ 
[\lambda]\!\in\!\Pro^{t-d(\cA)-2}_\R\setminus 
H_\cA\right\}}\!\subset\!\R^{t-d(\cA)-1}$, \hfill\mbox{}\\ 
according as $\cA$ is defective or not. We call any connected component $\cC$ 
of $\R^{t-d(\cA)-1}\setminus\Gamma_\sigma(\cA,B)$ a {\em signed reduced 
chamber}. We also call $\cC$ an {\em outer} or {\em inner chamber}, according 
as $\cC$ is unbounded or bounded. 
\dia 
\end{dfn}
\begin{ex} 
Continuing Example \ref{ex:penta}, there are $16$ possible 
choices for $\sigma$, if we identify\linebreak 
\scalebox{.96}[1]{sign sequences with their negatives. Among 
these choices, there are $11$ $\sigma$ yielding 
$\Gamma_\sigma(\cA,B)\!=\!\emptyset$.}\linebreak  
The remaining choices, along with their respective $\Gamma_\sigma(\cA,B)$ 
are drawn below. \dia 
\end{ex} 

\noindent 
\scalebox{.92}[.92]{
\begin{picture}(200,90)(12,-17)
\put(10,-15){
\fbox{\epsfig{file=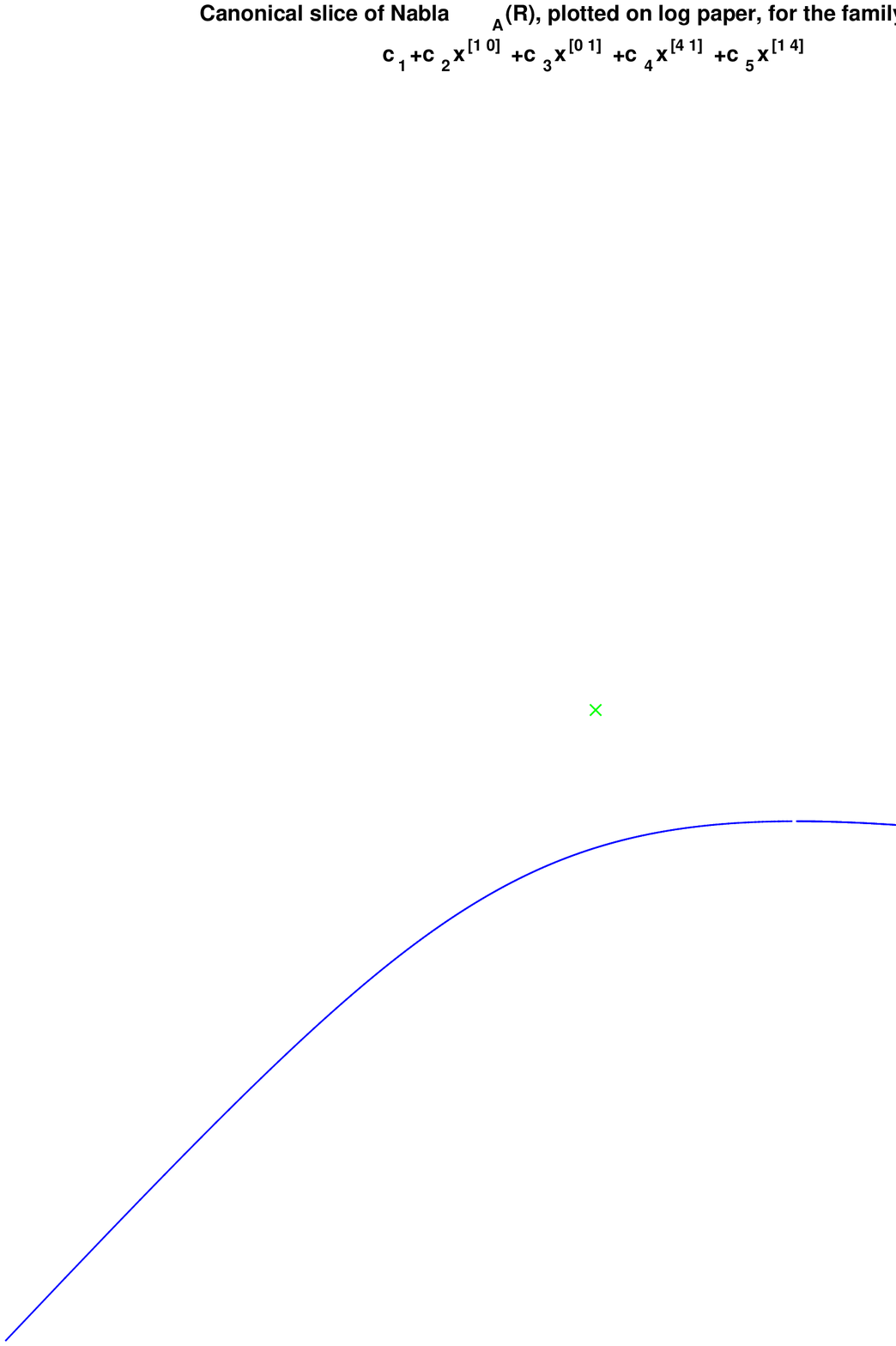,height=1.2in,clip=}}
\fbox{\epsfig{file=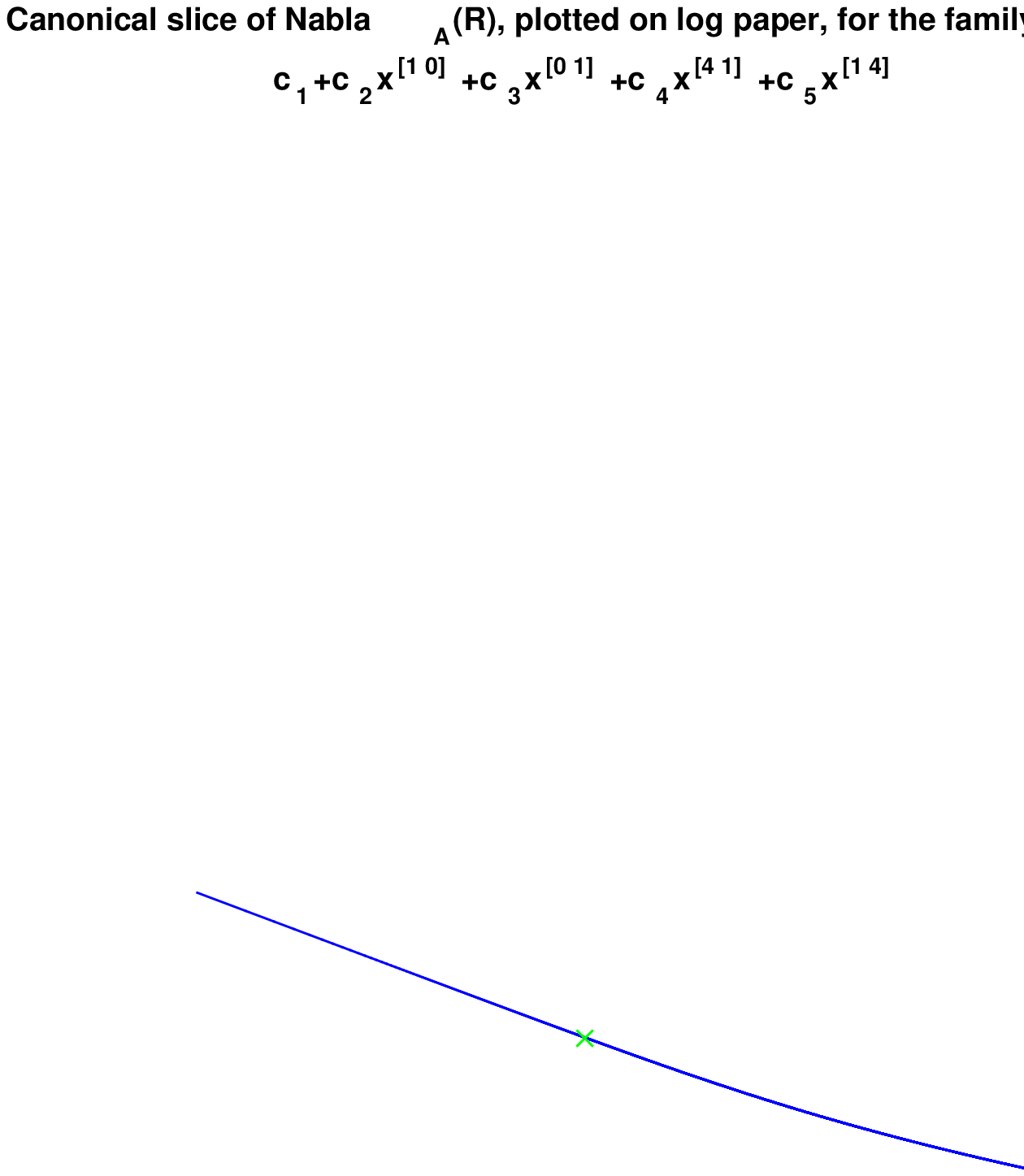,height=1.2in,clip=}}
\fbox{\epsfig{file=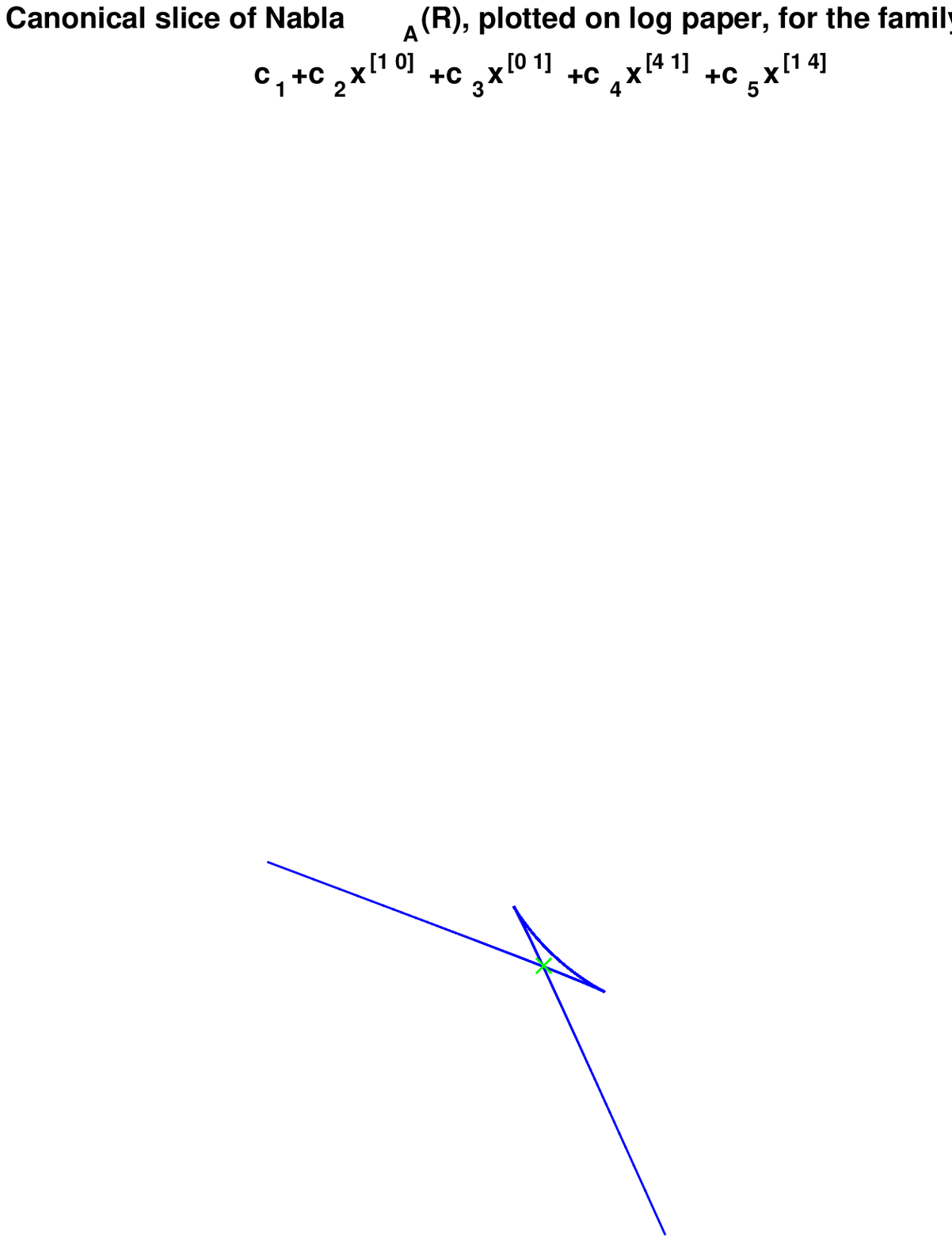,height=1.2in,clip=}}
\fbox{\epsfig{file=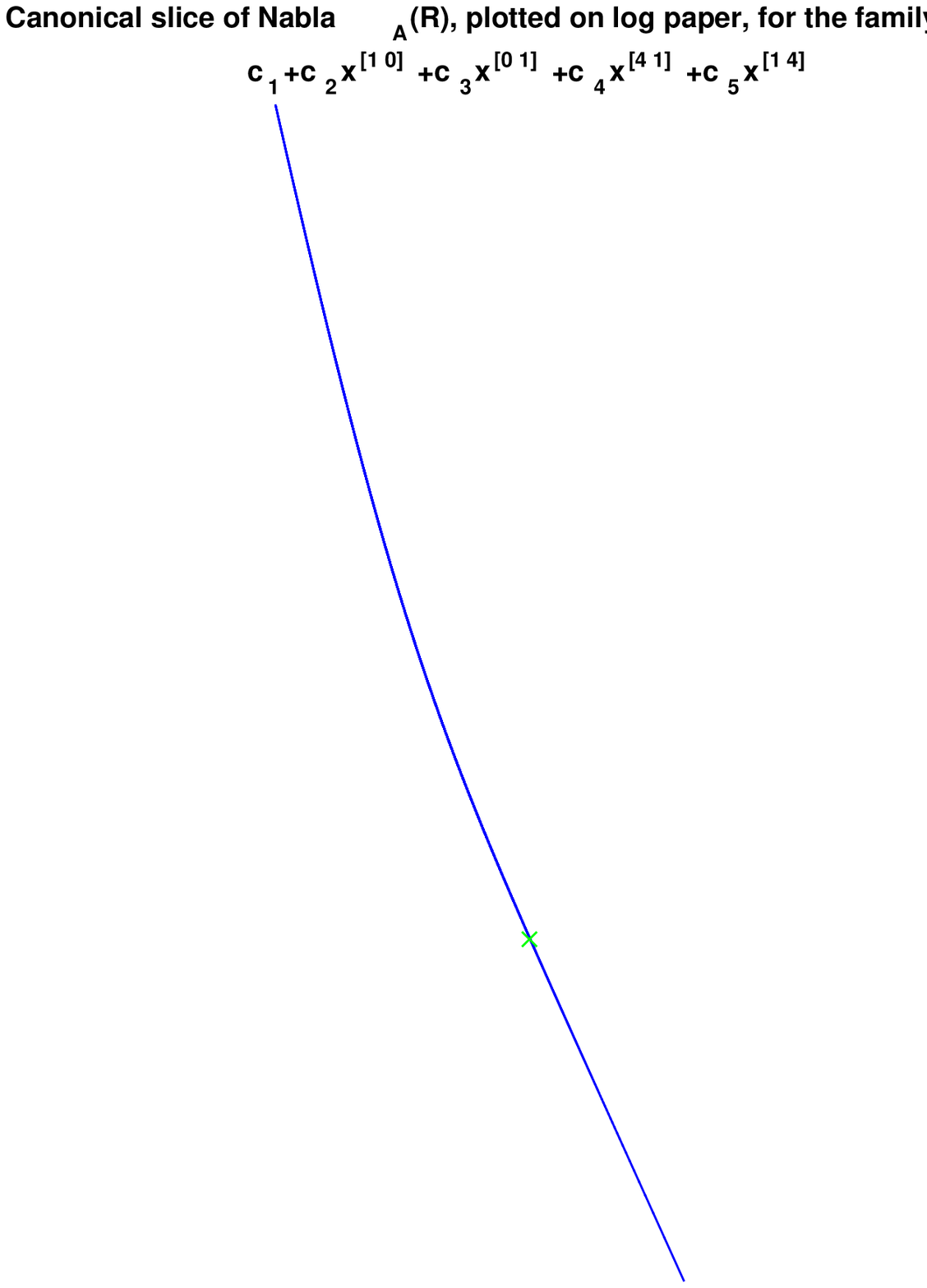,height=1.2in,clip=}}
\fbox{\epsfig{file=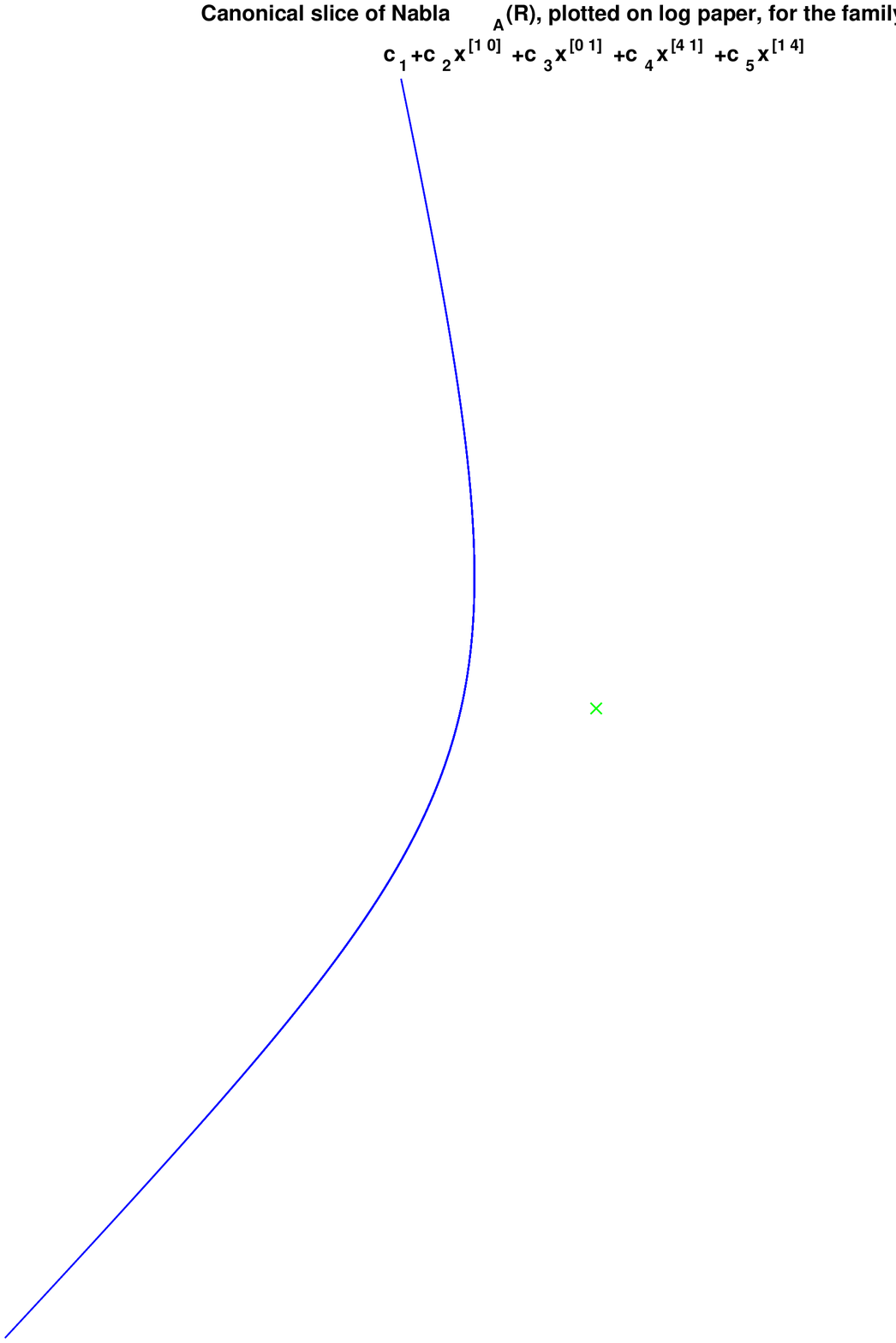,height=1.2in,clip=}}}
\put( 30,60){$++--+$}
\put(128,60){$-+++-$}
\put(230,-10){$+--++$}
\put(327,-10){$-++-+$}
\put(455,-10){$+-++-$}
\end{picture}}

\noindent 
Note that the curves drawn above are in fact unbounded, 
so the number of reduced signed chambers for the $\sigma$ above, from left to 
right, is respectively $2$, $2$, $3$, $2$, and $2$. (The tiny $\times$ in 
each illustration indicates the origin in $\R^2$.) In particular, 
only $\sigma\!=\!(1,-1,-1,1,1)$ yields an inner chamber. Note also that 
$\Gamma(\cA,B)$ is always the union of all the $\Gamma_\sigma(\cA,B)$.  
\begin{rem} While the shape of the reduced signed chambers 
certainly depends on the choice of $B$, the hyperplane arrangement 
$H_\cA$ and the number of signed chambers 
for any fixed $\sigma$ are independent of $B$. In particular, working 
with the $\Gamma_\sigma(\cA,B)$ (which have dimension $t-d(\cA)-2$) 
helps us visualize and work with the real part of $\Xi_\cA$ (which has 
dimension $t-2$ and involves up to $2^{t-1}$ orthants). \dia 
\end{rem} 

\noindent 
\begin{minipage}[t]{.6 \textwidth}
\vspace{0pt}
It is also important to note that innner chambers are where the isotopy type 
of $Z_\R(g)$ becomes more subtle, and no longer 
approachable via classical patchworking \cite{viropatch}. 
For instance, the $3$ possible isotopy types for $Z_\R(g)$ 
in our last example (with $\sigma\!=\!(1,-1,-1,1,1)$) are drawn  
within the boxes to the right. Note in particular that the isotopy type 
with a compact oval is {\em not} obtainable by applying patchworking 
to the point set $\cA$ in any obvious way. 

\mbox{}\hspace{.5cm}We close this section with the key reason we introduced 
$\Log|\cdot|$ earlier: it makes studying the curvature of the image of 
$\xi_{\cA,B}$ much easier. 
\end{minipage}  
\begin{minipage}[t]{.4 \textwidth}
\vspace{0pt}
\scalebox{.5}[.5]{\vbox{
\begin{picture}(300,380)(0,-60)
\put(0,0){\epsfig{file=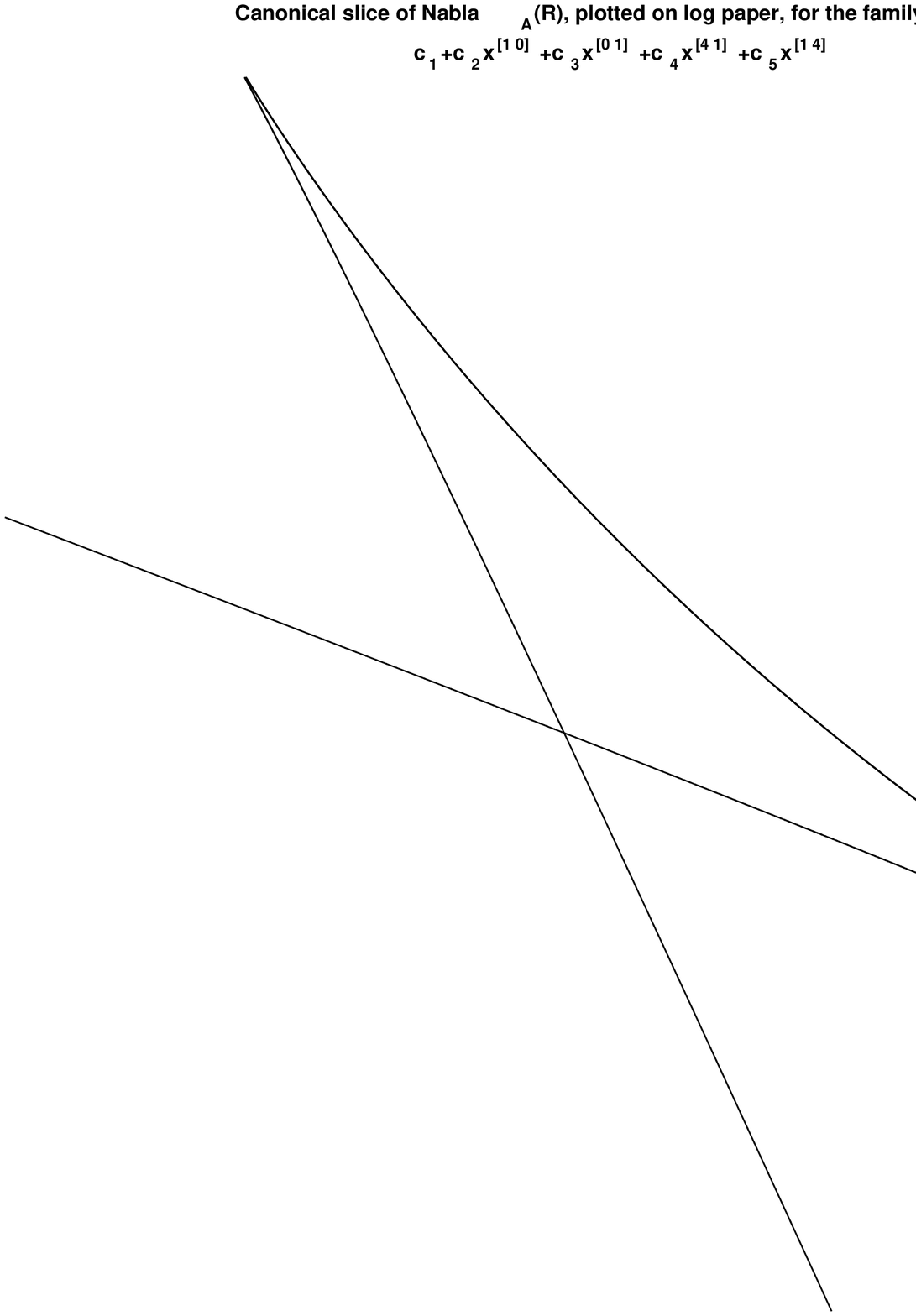,height=5in,clip=}}
\put(250,280){\fbox{\epsfig{file=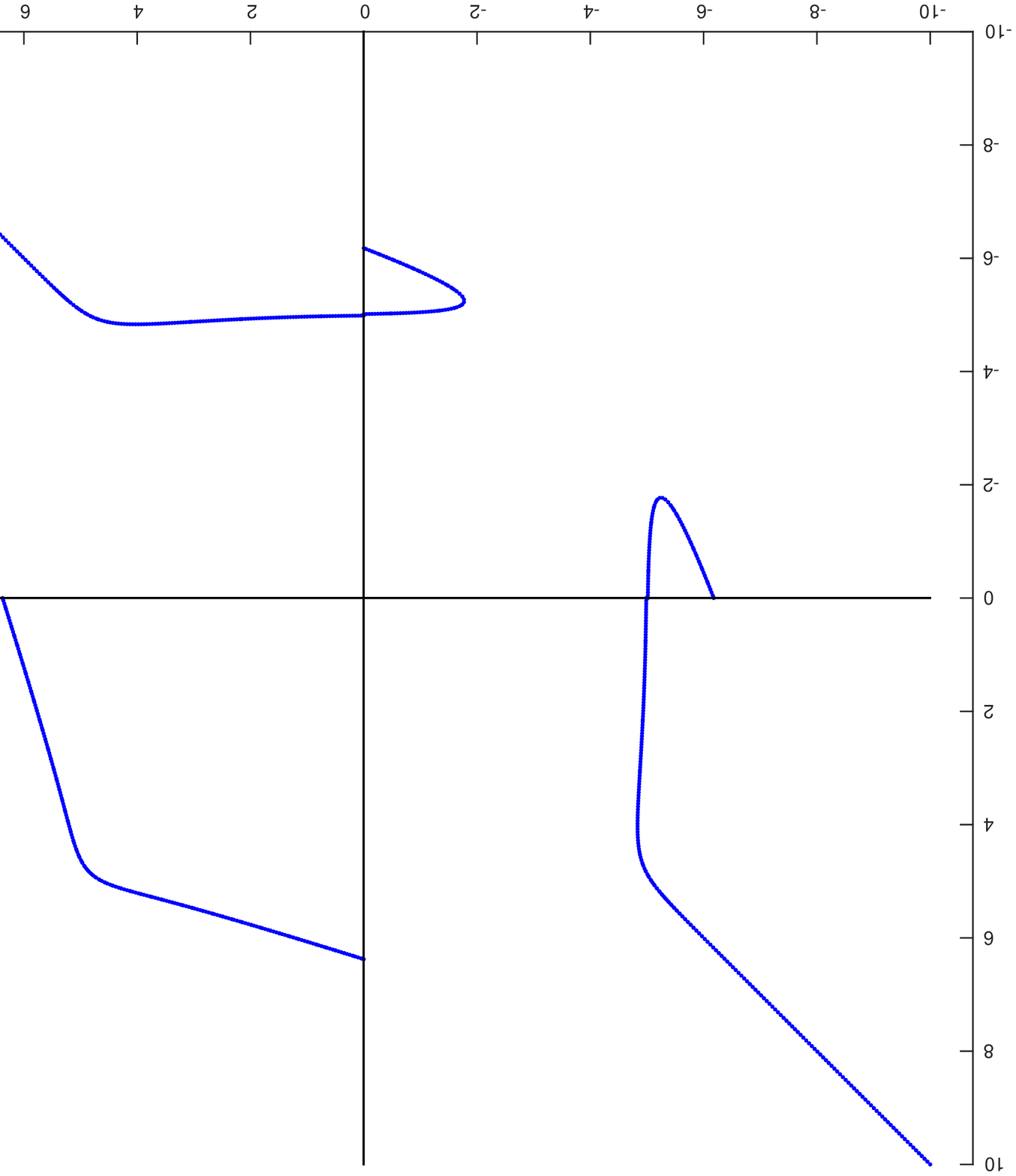,height=1.5in,angle=180,clip=}}}
\put(15,65){\fbox{\epsfig{file=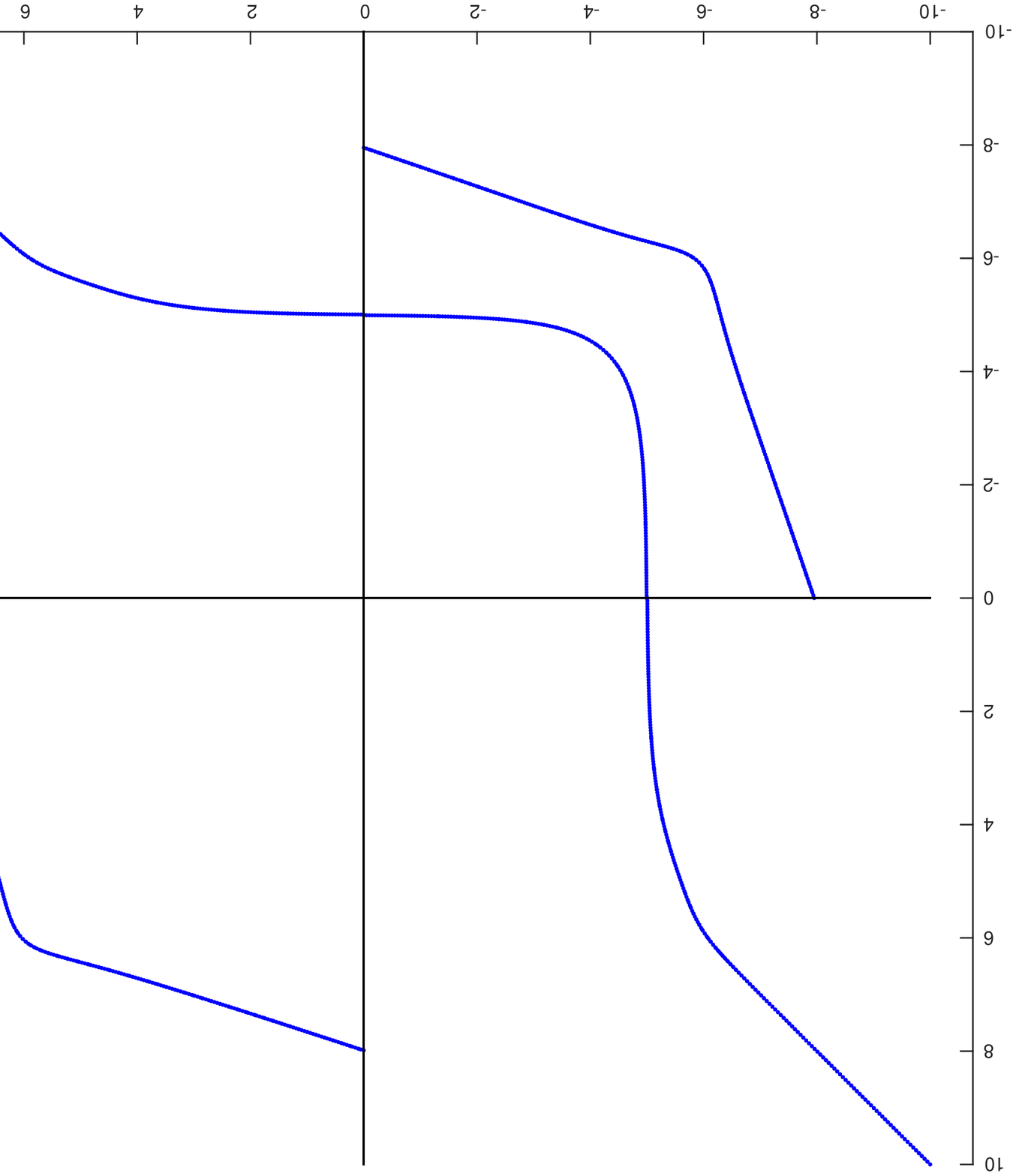,height=1in,angle=180,clip=}}}
\put(120,165){\fbox{\epsfig{file=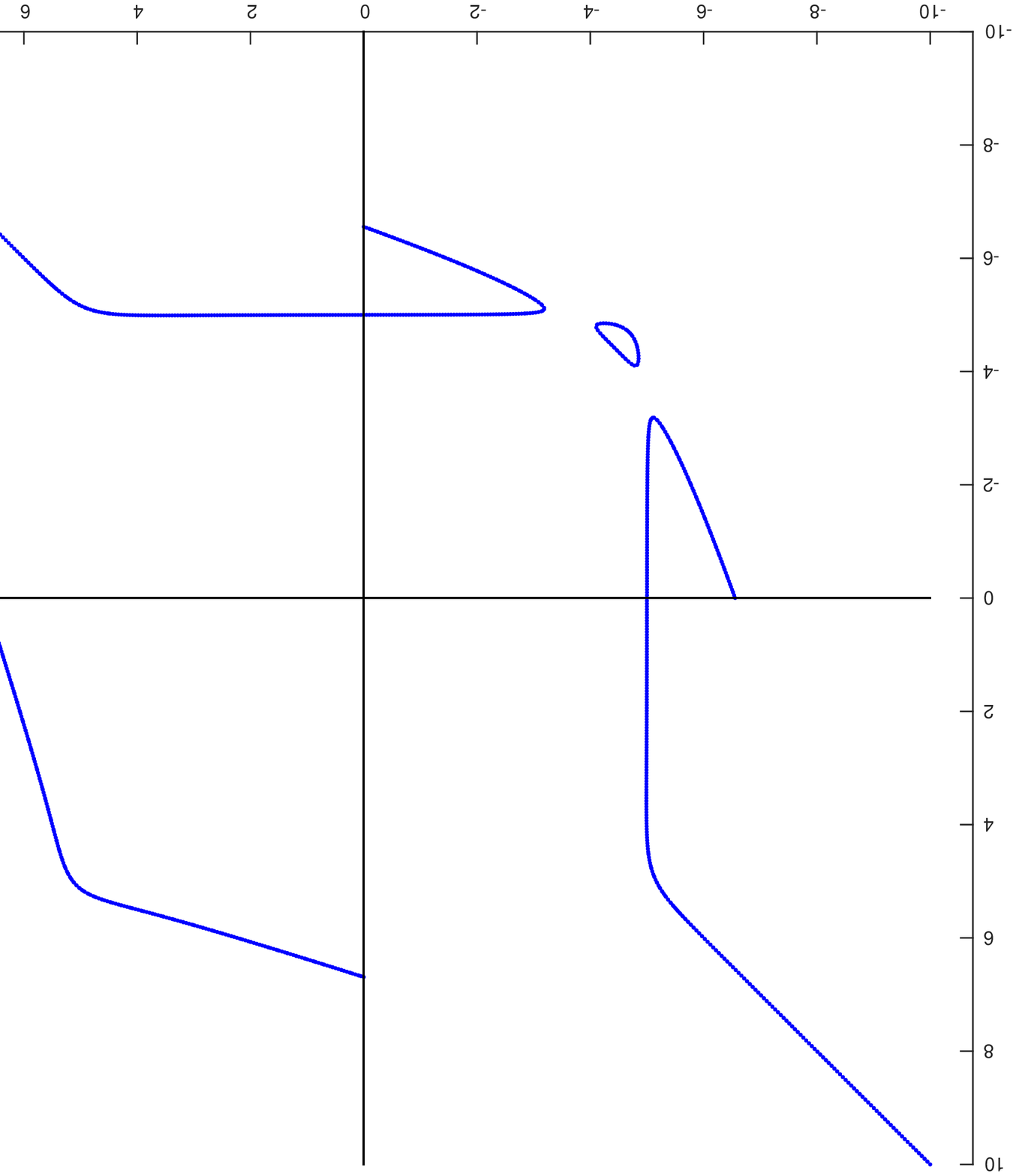,height=.8in,angle=180,clip=}}}
\end{picture}}} 
\end{minipage} 

\begin{thm} 
\label{thm:normal} 
Suppose $\cA\!\in\!\R^{n\times t}$ is non-defective. 
Then, at any point $[\ell]\!\in\!\Pro^{t-d(\cA)-2}_\R\setminus H_\cA$ where 
$\xi_{\cA,B}$ is differentiable, we have that 
$v\!\in\!\R^{t-d(\cA)-1}\setminus\{\bO\}$ is a normal vector to 
$\xi_{\cA,B}([\ell]) \Longleftrightarrow [v]\!=\![\ell]$. \qed 
\end{thm}

\noindent 
Recall that the {\em Gauss map} from $\Gamma(\cA,B)$ to $\Pro^1_\R$ 
is the map to $\Pro^1_\R$ obtained by mapping a point of 
$\Gamma(\cA,B)$ to the direction of its normal line. 
A simple consequence of Theorem \ref{thm:normal} is that the Gauss map 
in our setting is an injection with dense image. (An example is our 
illustration of $\Gamma(\cA,B)$ from Example \ref{ex:penta}.) 
The special case $\cA\!\in\!\Z^{n\times t}$ of Theorem \ref{thm:normal}
was already observed by Kapranov in \cite[Thm.\ 2.1 (b)]{kapranov}. The proof
of the more general Theorem \ref{thm:normal} is in fact almost
identically, since it ultimately reduces to elementary identities
involving linear combinations of logarithms in linear forms in $\lambda$, 
where one merely needs $\cA$ to be real.

\section{Morse Theory, Completed Contours, and Roots at Infinity} 
\label{sec:morse} 
Let us call $\cA\!\in\!\R^{n\times t}$ {\em combinatorially simplicial}
if and only if $\cA\cap Q$ has cardinality $1+\dim Q$ for every face $Q$ of
$\conv \{a_1,\ldots,a_t\}$.
(The books \cite{grunbaum,ziegler} are excellent standard references
on polytopes, their faces, and their normal vectors.)
Note that $\conv \{a_1,\ldots,a_t\}$
need {\em not} be a simplex for $\cA$ to be combinatorially simplicial
(consider, e.g., Example \ref{ex:penta}). We now state the main reason we care
about reduced signed chambers. In what follows, for $g(y)\!=\!\sum^t_{j=1} 
c_j e^{a_j\cdot y}$, we let $c_g\!:=\!(c_1,\ldots,c_t)$. 
\begin{thm}
\label{thm:morse} 
Suppose $\cA\!\in\!\R^{n\times t}$ is combinatorially simplicial, 
and $g_1$ and $g_2$ are each $n$-variate exponential
$t$-sums with spectrum $\cA$ and smooth real zero set. 
Suppose further that
$\sign(c_{g_1})\!=\!\pm \sign(c_{g_2})$, and $(\Log|c_{g_1}|)B$ and
$(\Log|c_{g_2}|)B$ lie in the same signed reduced discriminant chamber. Then
$Z_\R(g_1)$ and $Z_\R(g_2)$ are ambiently isotopic in $\Rn$. 
\end{thm}

\noindent 
We prove Theorem \ref{thm:morse} after recalling one definition and 
one lemma on a variant of the classical {\em moment map} \cite{smale,
souriau} from symplectic geometry. 
The special case $\cA\!\in\!\Z^{n\times t}$ of Theorem \ref{thm:morse}, 
without the use of $\Log$ or $B$, is alluded to near the beginning of
\cite[Ch.\ 11, Sec.\ 5]{gkz94}.
However, Theorem \ref{thm:morse} is really just an instance of
{\em Morse Theory} \cite{milnor,smt},
once one considers the manifolds defined by the fibers of the
map\linebreak $Z_\R(g)\mapsto (\Log|c_g|)B$ along paths inside a fixed
signed chamber. In particular, the assumption that $\cA$ be combinatorially
simplicial forces any topological change in $Z_\R(g)$ to arise solely from
singularities of $Z_\R(g)$ in $\Rn$. When $\cA$ is more general, topological
changes in $Z_\R(g)$ can arise from pieces of $Z_\R(g)$ approaching infinity,
with no singularity appearing in $\Rn$. So our chambers will eventually 
need to be cut into smaller pieces. 

So let us now make the notion of roots at infinity rigorous. 
\begin{dfn} 
\label{dfn:strat}
Let $\inte(S)$ denote the topological interior of any set $S\!\subseteq\!\Rn$, 
and let $\relint(Q)$ denote the relative interior of any $d$-dimensional
polytope $Q\!\subset\!\Rn$, i.e., $Q\!\setminus\!R$ where
$R$ is the union of all faces of $Q$ of dimension strictly
less than $d$ (using $\emptyset$ as the only face of
dimension $<\!0$). For any $w\!\in\!\Rn\setminus\{\bO\}$, 
we define $Q^w\!:=\!\{x\!\in\!Q \; | \; x\cdot w\!=\!\min_{y\in Q} 
y\cdot w\}$ to be the {\em face of $Q$ with inner normal $w$}. Finally, 
for any real exponential sum $g(y)\!=\!\sum^t_{j=1} c_j e^{a_j\cdot y}$, 
let $\init_w(g)\!:=\!\sum_{a_j\in \conv\{a_1,\ldots,a_t\}^w} c_je^{a_j 
\cdot y}$ be the {\em initial term summand of $g$ with respect to the 
weight $w$}.  \dia 
\end{dfn}  
\begin{lemma} \cite[Lemma 14]{tri}\footnote{The statement in 
\cite{tri} is for polynomials with real exponents, but is easily 
seen to be equivalent to the statement here via the substitution 
$x_i\!=\!e^{y_i}$.}   
\label{lemma:moment}
Given any $n$-dimensional convex compact polytope $P\!\subset\!\Rn$, there
is a real analytic diffeomorphism $\mu_P : \Rn \longrightarrow 
\inte(P)$. In particular, if $g$ is an $n$-variate exponential $t$-sum with 
spectrum having convex hull $P$ of dimension $n$, and
$w\!\in\!\Rn\!\setminus\!\{\bO\}$, then
$\mu_P(Z_\R(g))$ has a limit point in $\relint(P^w) \Longrightarrow 
\init_w(g)$ has a root in $\Rn$. Moreover, there is a
real analytic diffeomorphism between $Z_\R(\init_w(g))\!\subset\!\Rn$ and
$\left(\relint(P^w)\cap\overline{\mu_P(Z_\R(g))}\right)
\times \R^{n-\dim P^w}$. \qed 
\end{lemma}

\noindent 
Note that the converse of Lemma \ref{lemma:moment} need not hold: A simple
counter-example is $g(y)\!:=\!(e^{2y_1}+e^{2y_2}-1)^2+(e^{y_1}-1)^2$ 
and $w\!=\!(0,1)$. In what follows, recall that $n$-manifolds (resp.\ 
$n$-manifolds with boundary) are defined by coordinate charts that are 
diffeomorphic to a (relatively) open subset of $\Rn$ (resp.\ $\R^{n-1}
\times (\R_+\cup\{0\})$) (see, e.g., \cite{hirsch}). More generally, an 
{\em $n$-manifold with corners} is defined via charts in $(\R_+\cup\{0\})^n$ 
instead (see, e.g., \cite{smt}). Manifolds with corners admit a natural 
Whitney stratification into (smooth, open) sub-manifolds without boundary. 
In fact, a by-product of Definition \ref{dfn:strat} is just such a 
stratification via the face lattice of $P$.

\medskip 
\noindent 
{\bf Proof of Theorem \ref{thm:morse}:} Let $P\!:=\!\conv\{a_1,
\ldots,a_t\}$, let $\cC$ be the unique signed reduced 
chamber containing $(\Log|c_{g_1}|)B$ and $(\Log|c_{g_2}|)B$, 
let $\bC$ be the fiber of $\Log|\cdot||_{\Pro^{t-1}_\sigma}$ 
over $\cC$, let $\fii\!=\![\fii_1:\cdots :\fii_t] 
: [0,1] \longrightarrow \bC$ be any analytic 
path connecting $[c_{g_1}]$ and $[c_{g_2}]$, and let\\  
\mbox{}\hfill $M\!:=\!\left\{(y,t)\!\in\!\R^n\times [0,1]\; \left| \; 
\sum^t_{j=1}\fii_j(t) e^{a_j\cdot y}\!=\!0\right.\right\}$.\hfill\mbox{}\\ 
Also let $p : M\longrightarrow [0,1]$ 
be the natural projection that forgets the first $n$ coordinates 
of $M$. It is then easily checked that $p$ has no critical points. 
Furthermore, letting $\bM$ denote the Euclidean closure of $\mu_P(M)$ 
in $P$, it easily checked via Lemmata \ref{lemma:basic} and 
\ref{lemma:moment} that $\bM\cap Q$ is a 
smooth hypersurface in $Q$, for any face $Q$ of $P$. So then, $\bM$ is a 
compact manifold with corners, and $p$ extends naturally to a 
projection $\bar{p} : \bM\longrightarrow [0,1]$ that also has no critical 
points. So then, by the Regular Interval Theorem (see, e.g., 
\cite[Thm.\ 2.2, Pg.\ 153]{hirsch} for the case of manifolds with boundary), 
$\bM$, $\overline{\mu_P(Z_\R(g_1))}
\times [0,1]$, and $\overline{\mu_P(Z_\R(g_2))}\times[0,1]$ are 
diffeomorphic to each other. In other words, the ends of $\bM$ 
($\mu_P(Z_\R(g_1))$ and $\mu_P(Z_\R(g_2))$) are ambiently isotopic in $P$, 
and thus $Z_\R(g_1)$ and $Z_\R(g_2)$ are ambiently isotopic in $\Rn$. \qed 

\medskip  
To address arbitrary $\cA$ we'll first need a little more terminology.
\begin{dfn}
\label{dfn:chamber2}
Given any $\cA\!\in\!\R^{n\times t}$ with distinct
columns $a_1,\ldots,a_t$, and any inner normal $w\!\in\!\Rn$ to a face of
$\conv \{a_1,\ldots,a_t\}$, we let $\cA^w\!:=\![a_{j_1},\ldots,a_{j_r}]$ 
denote the sub-matrix of $\cA$ corresponding to the set
$\{a\!\in\!\cA \; | \; a\cdot w\!=\!\min_{a'\in\cA}\{a'\cdot w\}\}$.
We call $\cA^w$ a {\em (proper) non-simplicial face} of $\cA$ when
$d(\cA^w)\!\leq\!d(\cA)-1$ and $\cA^w$ has at least $d(\cA^w)+1$ columns.
Also let $B^w$ be any matrix whose columns form a basis for the right
nullspace of $\widehat{\left(A^w\right)}$, and let $\pi_w : \C^t  
\longrightarrow \C^r$ be the natural coordinate projection map defined by
$\pi_w(c_1,\ldots,c_t)\!:=\!(c_{j_1},\ldots,c_{j_r})$. When $\cA$ is
non-defective we then define the {\em completed} reduced signed contour,\\ 
\mbox{}\hfill
$\widetilde{\Gamma}_\sigma(\cA,B)\subset\R^{t-d(\cA)-1}$,\hfill\mbox{}\\ 
to be the union of $\Gamma_\sigma(\cA,B)$ and\\
\mbox{}\hfill $\displaystyle{\bigcup\limits_{\substack{\cA^w \text{a non-}\\ 
\text{simplicial}\\ \text{face of } \cA}} 
\overline{\left\{\left.\pi^{-1}_w\!\left(\Log|\lambda(B^w)^\top|\right)B 
\; \right| \; \mathrm{sign}\left(\lambda (B^w)^\top\right) 
\!=\!\pm \pi_w(\sigma)  \ , \ [\lambda]\!\in\!
\Pro^{t-d(\cA)-2}_\R\setminus H_\cA\right\}}}$.\hfill\mbox{}\\
We call any unbounded connected component of $\R^{t-d(\cA)-1}\setminus
\widetilde{\Gamma}_\sigma(\cA,B)$ a {\em refined} outer chamber. 
Finally, we call $\widetilde{\Gamma}(\cA,B)\!:=\!\bigcup\limits_{\sigma\in
\{\pm 1\}^{t}}\widetilde{\Gamma}_\sigma(\cA,B)$ a {\em completed 
reduced (unsigned) contour}. \dia
\end{dfn}

\begin{ex}
\label{ex:inf}
When $\cA\!=$\scalebox{.6}[.6]{$\begin{bmatrix}0 & 1 & 0 & 2 & 0\\ 
0 & 0 & 1 & 0 & 2 \end{bmatrix}$} it is easy to find a $B$ yielding the
following reduced contour $\Gamma(\cA,B)$ and completed reduced contour
$\widetilde{\Gamma}(\cA,B)$:\\
\mbox{}\hfill\epsfig{file=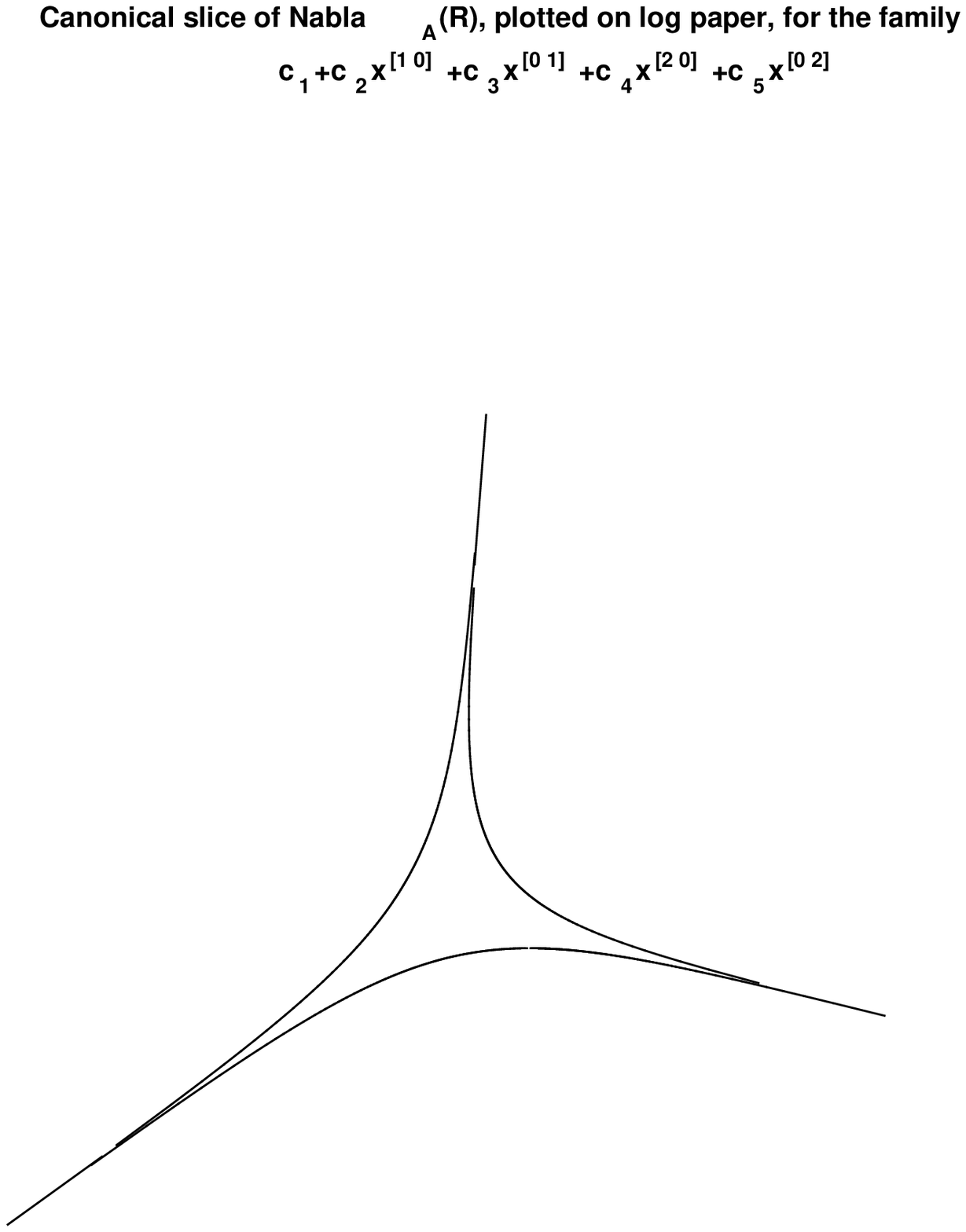,height=1.6in,clip=}
\hfill\epsfig{file=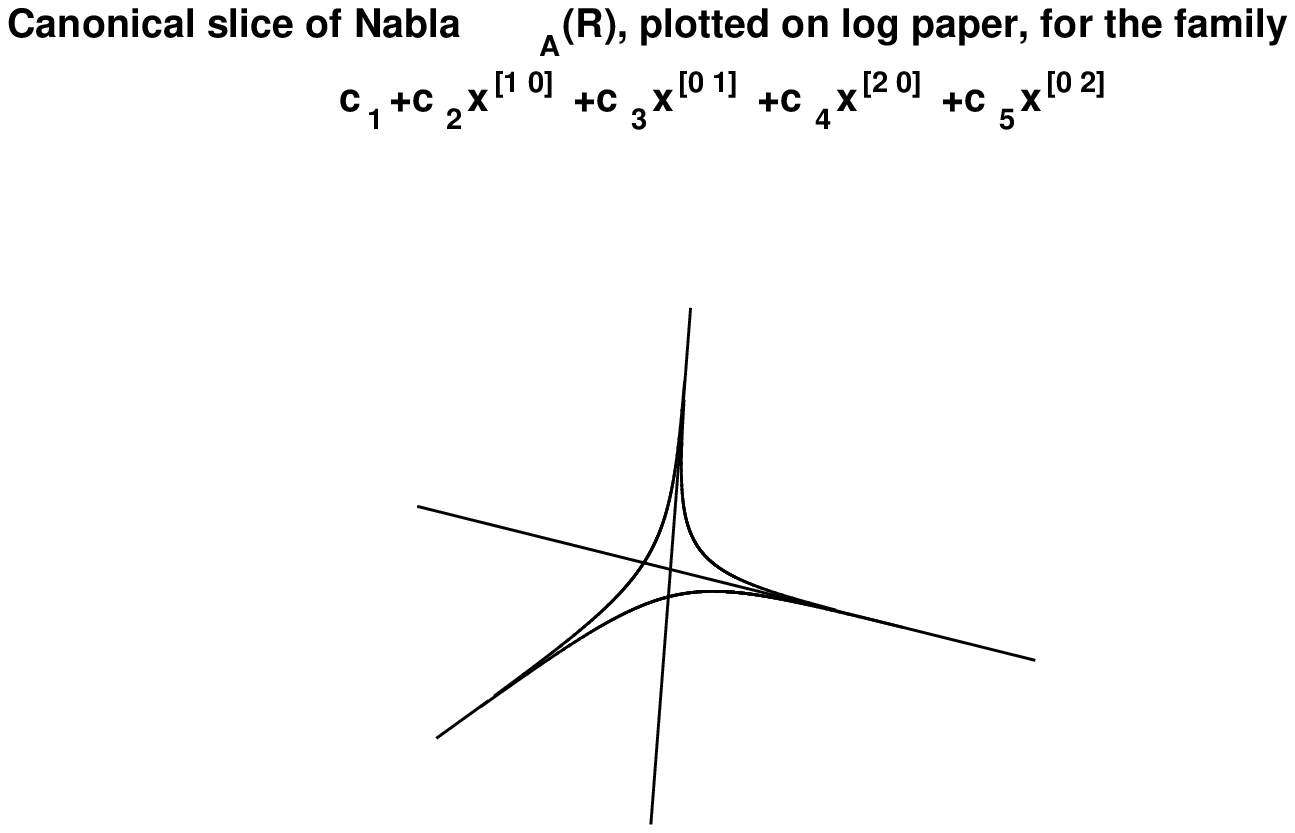,height=1.6in,clip=}
\hfill\mbox{}\\
Note in particular that $\widetilde{\Gamma}(\cA,B)\!=\!\Gamma(\cA,B)\cup 
S_1\cup S_2$ where $S_1$ and $S_2$ are lines that can be viewed as
line bundles over points. These points are in fact $(\Log|\Xi_{\cA_1}|)B$ and
$(\Log|\Xi_{\cA_2}|)B$ where $\cA_1$ and $\cA_2$ are the facets of $\cA$ with
respective outer normals $(-1,0)$ and $(0,-1)$, and
$B\!\approx$\scalebox{.6}[.6]
{$\begin{bmatrix}[r] 0.4335 & -0.8035 & -0.0635 & 0.4018 & 0.0317\\ 
0.3127 & 0.2002 & -0.8256 & -0.1001 & 0.4128\end{bmatrix}^\top$}. \dia
\end{ex}
\begin{prop}
\label{prop:non}
If $\cA$ is non-defective and not combinatorially simplicial then
$\cA$ has at most $t-d(\cA)-1$ non-simplicial faces. \qed
\end{prop}
\begin{prop}
\label{prop:2face}
\scalebox{.95}[1]{Suppose $k\!=\!3$, $\cA$ has exactly $2$ non-simplicial
facets, $d(\cA)\!=\!n$, $B\!\in\!\R^{(n+3)\times 2}$}\linebreak
is any matrix whose columns form a basis for the right
nullspace of $\hA$, and $[\beta_{i,1},\beta_{i,2}]$ is\linebreak
the $i\thth$ row of $B$.
Then $\{[\beta_{i,1}:\beta_{i,2}]\}_{i\in\{1,\ldots,n+3\}}$ has cardinality
$n+1$ as a subset of $\Pro^1_\R$, and $\widetilde{\Gamma}(\cA,B)
\setminus\Gamma(\cA,B)$ is a union of $2$ lines. \qed
\end{prop}
\begin{thm} 
\label{thm:morse2}
Suppose $\cA\!\in\!\R^{n\times t}$, and $g_1$ and $g_2$ are each 
$n$-variate exponential $t$-sums with spectrum $\cA$ and smooth real zero 
set. Suppose further that 
$\sigma\!:=\!\sign(c_{g_1})\!=\!\pm \sign(c_{g_2})$,
and $(\Log|c_{g_1}|)B$ and $(\Log|c_{g_2}|)B$
lie in the same connected component of
$\R^{t-d(\cA)-1}\setminus\widetilde{\Gamma}_\sigma(\cA,B)$. Then
$Z_\R(g_1)$ and $Z_\R(g_2)$ are ambiently isotopic in $\Rn$. 
\end{thm}
\begin{ex}
Observe that the circle defined by
$\left(u+\frac{1}{2}\right)^2+\left(v-2\right)^2\!=\!1$ intersects the
positive orthant, while the circle defined by
$\left(u+\frac{3}{2}\right)^2+\left(v-\frac{3}{2}\right)^2=1$ does not.
Consider then $\cA\!=$\scalebox{.6}[.6]{$\begin{bmatrix}0 & 1 & 0 & 2 & 0\\ 
0 & 0 & 1 & 0 & 2 \end{bmatrix}$} as in our last example, and let
\textcolor{blue}{$g_1\!=\!\left(e^{y_1}+\frac{1}{2}\right)^2
+\left(e^{y_2}-2\right)^2-1$} and\linebreak
\textcolor{red}{$g_2\!=\!\left(e^{y_1}+\frac{3}{2}\right)^2+\left(e^{y_2}
-\frac{3}{2}\right)^2-1$}. Then \textcolor{blue}{$g_1$} and \textcolor{red}
{$g_2$} have spectrum $\cA$, $\sign(g_1)\!=\!\sign(g_2)\!=\!\sigma$ with
$\sigma\!=\!(1,1,-1,1,1)$, and \textcolor{blue}{$(\Log|c_{g_1}|)B$} and
\textcolor{red}{$(\Log|c_{g_2}|)B$} lie in the same reduced signed
$\cA$-discriminant chamber (since
$\Gamma_\sigma(\cA,B)\!=\!\emptyset$ here). However, \textcolor{blue}
{$Z_\R(g_1)$} consists of a single smooth arc, while
\textcolor{red}{$Z_\R(g_2)$} is empty.
This is easily explained by the {\em completed} contour
$\widetilde{\Gamma}_\sigma(\cA,B)$ consisting of two lines, and
\textcolor{blue}{$(\Log|c_{g_1}|)B$} and \textcolor{red}{$(\Log|c_{g_2}|)B$}
lying in distinct {\em refined} chambers  
as shown, respectively via the symbols $\circ$ and $*$, below to
the right. \dia
\end{ex}

\vspace{-.2cm}
\noindent
\begin{minipage}[t]{.7 \textwidth}
\vspace{0pt}
{\bf Proof of Theorem \ref{thm:morse2}:} We follow exactly the same 
set-up as the proof of Theorem \ref{thm:morse}. However, although 
$\cA$ may not be combinatorially simplicial, our use of 
$\widetilde{\Gamma}_\sigma(\cA,B)$ instead of $\Gamma_\sigma(\cA,B)$ 
guarantees that $\bar{p}$ is still a projection from $\bM$ to $[0,1]$ 
without critical points. So we can again obtain an ambient isotopy via the 
Regular Interval Theorem as before. \qed 
\end{minipage}\hspace{1cm}
\begin{minipage}[t]{.2 \textwidth}
\vspace{0pt}
\begin{picture}(0,0)(0,0)
\put(-25,-60){\epsfig{file=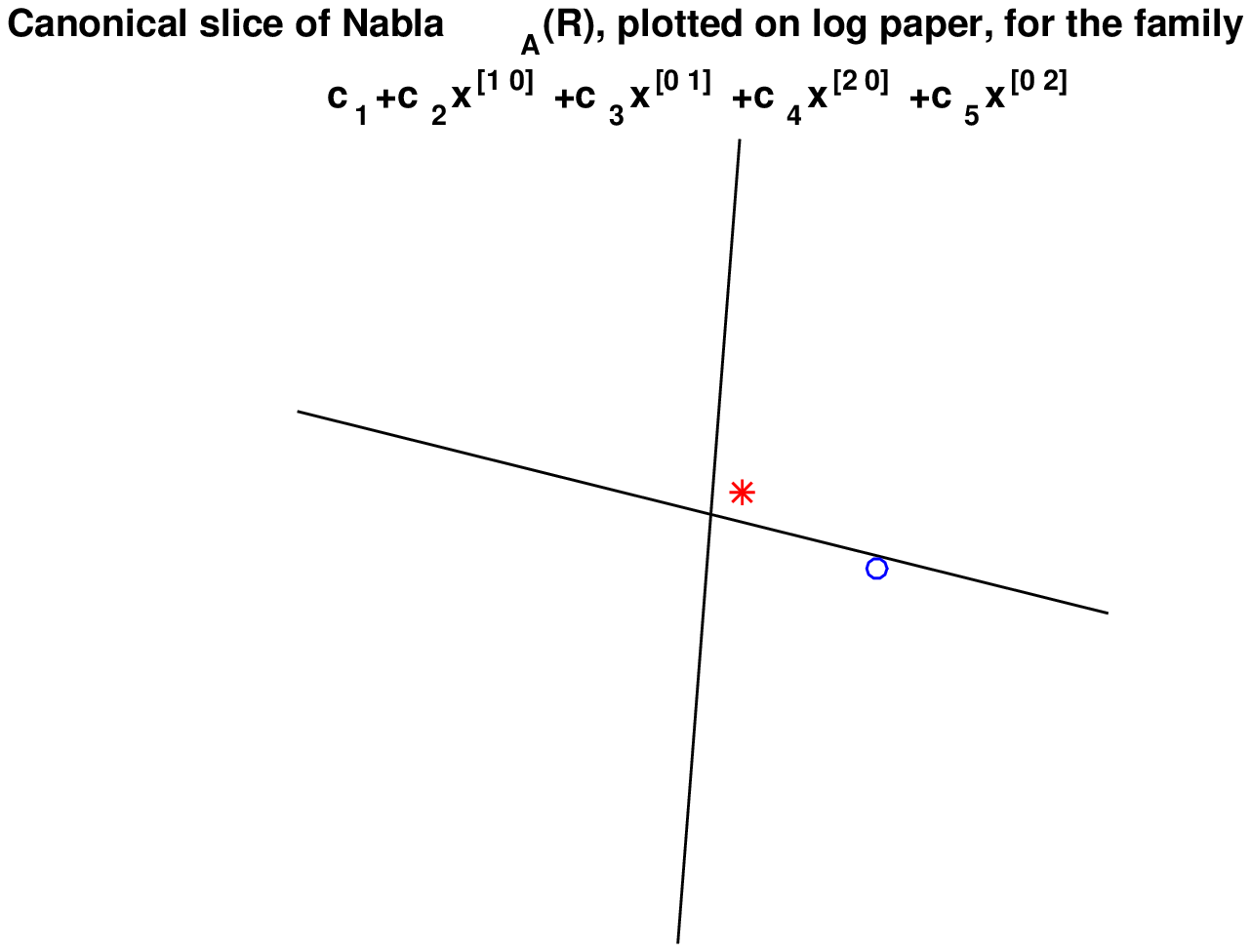,height=.9in,clip=}}
\end{picture}
\end{minipage}

\begin{rem} An alternative to our use of the moment map would have been 
to instead embed $\Rn$ into a suitable (compact) {\em quasi-fold} (see, e.g., 
\cite{prato,nonrat}). Roughly, quasi-folds are to toric varieties what 
real exponents are to integral exponents. \dia 
\end{rem} 

We are now ready to prove our main results. 

\section{The Proof of Our Parametrization: Theorem \ref{thm:gena}} 
\label{sec:back} 
Let $c\!:=\!(c_1,\ldots,c_t)\!\in\!(\C\setminus\{0\})^t$. 
Rewriting the equations defining $g(y)\!=\!\frac{dg}{dy_1}\!=\cdots=
\!\frac{dg}{dy_n}\!=\!0$, we see that $g$ has a singularity at $y\!\in\!\Cn$ 
if and only if\\ 
\mbox{}\hfill $\hA\begin{bmatrix}c_1e^{a_1\cdot y}\\ \vdots \\ c_te^{a_t 
\cdot y} \end{bmatrix} \!=\!\bO$.\hfill \mbox{}\\ 
The last equality is equivalent to $c\odot e^{yA}\!=\!\lambda 
B^\top$ for some $[\lambda]\!\in\!\left.\Pro^{t-d(\cA)-1}_\C\right\backslash 
H_\cA$, thanks to the definition of the right nullspace of $\hA$. 
(Note that Proposition \ref{prop:pyr} implies that 
$\lambda$ must be such that no coordinate of $\lambda B^\top$ 
vanishes.) Dividing the $i\thth$ coordinate 
of each side by\linebreak  
\scalebox{.97}[1]{$e^{a_j\cdot y}$, we thus obtain that the image 
$\psi_\cA\!\left(\left(\left.\Pro^{t-d(\cA)-2}_\C\right\backslash 
\cH_\cA\right)\times \Cn\right)$ is exactly 
$\Xi_\cA\setminus\{c_1\cdots c_t\!=\!0\}$.}\linebreak 
So the last two sets have the same Euclidean closure.
Note in particular that $\psi_\cA$ is analytic on 
$\left(\left.\Pro^{t-d(\cA)-2}_\C\right\backslash \cH_\cA\right)
\times \Cn$, and $\left(\left.\Pro^{t-d(\cA)-2}_\C\right\backslash  
\cH_\cA\right)\times \Cn$ is connected. (The latter assertion follows from 
the fact that complements 
of algebraic hypersurfaces in $\C^N$ are connected: See, e.g., 
\cite[Pg.\ 196]{bcss}.) So 
$\psi_\cA\!\left(\left(\left.\Pro^{t-d(\cA)-2}_\C\right\backslash 
\cH_\cA\right)\times \Cn\right)$ is a connected analytic hypersurface 
and we thus obtain Assertion (1).   

Let us now assume that all the entries of $\cA$ are real. Since 
$\hA B\!=\!\bO$, we 
immediately obtain that $\Log\!\left|\Xi_\cA\cap \Pro^{t-1}_\sigma\right|$ 
is an $\Rn$-bundle over $Z\!:=\!\left(\Log\left|\Xi_\cA\cap\Pro^{t-1}_\sigma
\right|\right)B$, with $\Gamma_\sigma(\cA)$ 
and $Y_\sigma$ subsets of $Z$ to be described shortly. 
So we now prove Assertions (a) and (b). 

Let $\Sigma\!\subset\!\Pro^{t-d(\cA)-2}_\R$ denote the singular locus 
of $\xi_{\cA,B}$. Since $\xi_{\cA,B}([\lambda],y)\!=\!
(\Log|\psi_\cA([\lambda,y])|)B$ for 
all $\lambda$ and $y$, Theorem \ref{thm:normal} implies that 
$\Gamma_\sigma(\cA,B)$ is a dense open subset of $Z$, provided 
we prove that $\Sigma$ is an algebraic set of positive 
codimension. This is indeed the case: By an elementary calculation, 
$\xi_{\cA,B}$ is singular at $[\lambda]\!\in\!\Pro^{t-d(\cA)-1}_\C$ 
if and only if some $(t-d(\cA)-2)\times (t-d(\cA)-1)$ sub-matrix of 
the $(t-d(\cA)-1)\times (t-d(\cA)-1)$ matrix 
\[ M(\lambda):=\begin{bmatrix} 
\left(b^\top_1\odot \left(\frac{1}{\lambda\cdot \beta^\top_1},
\ldots,\frac{1}{\lambda\cdot \beta^\top_t}\right)\right)B\\
\vdots \\ 
\left(b^\top_{t-d(\cA)-1}\odot \left(\frac{1}{\lambda\cdot \beta^\top_1},
\ldots,\frac{1}{\lambda\cdot \beta^\top_t}\right)\right)B\\
\end{bmatrix} 
\] 
has rank $<\!t-d(\cA)-2$, and this is clearly an algebraic condition 
(defined over $\R$ in fact). 
Letting $\Delta\!\subset\!\Pro^{t-d(\cA)-1}_\R$ be the open polyhedron 
consisting of all $\lambda$ such that $\sign(\lambda B^\top)\!=\!\pm \sigma$ 
we then immediately obtain Assertion (a) (with $\Gamma_\sigma(\cA)\!=\!
\xi_{\cA,B}(\Delta\setminus\Sigma)$), save for the 
stated quantitative bound: $\Delta\setminus\Sigma$ is the complement of 
a real algebraic hypersurface and thus has only finitely many connected 
components, and $\xi_{\cA,B}$ is analytic on these connected components. 
We also obtain Assertion (b) (with 
$Y_\sigma$ the codimension $\geq\!2$ part of $Z$), 
since any analytic set is a countable locally finite union of 
manifolds \cite{loja}. So we now only need to prove the quantitative bound 
from Assertion (a). 

Let $R(\lambda)\!:=\!(\det M(\lambda))\left(
\prod^t_{i=1} \beta_i\cdot \lambda\right)^{t-d(\cA)-1}$. 
Clearly, $\Sigma\!\subseteq\!Z_\R(R)$. Also, since 
$\det M(\lambda)$ consists of a linear 
combination of products of terms of the form  
$\sum^t_{i=1}\frac{\gamma_i}{\beta_i\cdot \lambda}$, $R$ is clearly homogeneous 
and has degree no greater than $(t-d(\cA)-1)(t-1)$. 
By the classical Oleinik-Petrovsky/Milnor-Thom bound
(see, e.g., \cite[Ch.\ 16, Pg.\ 308, Prop.\ 3]{bcss}), 
$\Delta\setminus\Sigma$ has no more than $(t-d(\cA)-1)(t-1)(2(t-d(\cA)-1)
(t-1)-1)^{2t-d-2}$ connected components. So we are done. \qed 

\section{Counting Isotopy Types: Proving Theorem \ref{thm:nearckt} 
and Lemma \ref{lemma:weird}} 
\label{sec:topo} 
We will first need the following convenient characterization of the 
singularities of a reduced contour when $t\!=\!n+3$. 
\begin{lemma} 
\label{lemma:cusp} 
If $\cA\!\in\!\R^{n\times (n+3)}$ is non-defective then any reduced 
$\cA$-discriminant contour has no more than $n$ cusps. 
\end{lemma} 

\noindent 
{\bf Proof:} Considering $\psi_\cA([\lambda_1:\lambda_2])\!=\!
(\psi_1([\lambda_1:\lambda_2]),\psi_1([\lambda_1:\lambda_2]))$ 
and dehomogenizing by setting $(\lambda_1,\lambda_2)\!=\!(1,\lambda)$, 
we can detect cusps by setting $\frac{\partial \psi_1}{\partial \lambda}
\!=\!\frac{\partial \psi_2}{\partial \lambda}\!=\!0$. 
One then obtains a pair of 
equations of the form $\sum^{n+3}_{i=1}\frac{r_i}{\beta_{i,1}
+\beta_{i,2}\lambda}\!=\!0$. Clearing denominators, we then 
obtain a polynomial of degree $\leq\!n+2$. One can then check that 
the two leading coefficients are exactly $0$, and thus we in fact 
obtain a polynomial of degree $\leq\!n$. So we clearly obtain that 
no more than $n$ points of the form $[\lambda_1:\lambda_2]\!\in\!\Pro^1_\R$ 
can yield a cusp. \qed  

\medskip 
In 1826, Steiner studied line arrangments in $\R^2$ and proved that 
$m$ lines determine no more than $\frac{m(m-1)}{2}+m+1$ connected 
components for the complement of their union \cite{steiner}. Thanks to 
our development so far, our signed reduced contour consists of a union of 
a single (possibly singular) arc $C$. So we will need an analogue of 
Steiner's count for this non-linear arrangement. Fortunately, Theorem 
\ref{thm:normal} and Lemma \ref{lemma:cusp} 
imply that we can easily derive such an analogue. In particular, 
Theorem \ref{thm:normal} implies that $C$ is locally convex, and 
its subarcs between cusps can intersect at most once (and only 
when the subarcs do not share a cusp). 
Employing a result from the second author's Ph.D.\ thesis, 
we then obtain the following: 
\begin{thm} 
\label{thm:korben} 
\cite[Thm.\ 3.6 \& 3.7]{rusek} 
Suppose $C\subset\!\R^2$ is a piece-wise smooth 
arc with exactly $\ell$ cusps. Suppose also that no two distinct smooth 
sub-arcs of $C$ intersect more than once, and that 
sub-arcs sharing a cusp do not intersect. 
Then $C$ has at most $\frac{\ell(\ell+1)}{2}
-(\ell+1)$ intersections, and the complement $\R^2\setminus C 
$ has at most $\frac{\ell(\ell+1)}{2}-\ell+1$ connected components. \qed 
\end{thm} 

Combining Lemma \ref{lemma:cusp} and Theorem \ref{thm:korben} we then 
obtain the following: 
\begin{cor} 
\label{cor:chamber} 
For any non-defective $\cA\!\in\!\R^{n\times (n+3)}$,  
the complement of the reduced $\cA$-discriminant contour consists of 
no more than $\frac{1}{2}n^2-\frac{1}{2}n+1$ connected components. \qed 
\end{cor} 

The final fact we'll need before proving Theorem \ref{thm:nearckt} is 
a refinement of Theorem \ref{thm:nearckt} involving 
$\cA\!\in\!\R^{n\times (n+2)}$: 
\begin{lemma} 
\label{lemma:ckt} 
Assume $\{a_1,\ldots,a_{n+2}\}\!\subset\!\Rn$
does not lie in any affine hyperplane and assume $\cA\!=\![a_1,\ldots,a_{n+2}]$ 
and $B\!\in\!\R^{(n+2)\times 2}$ satisfy the hypotheses of Definition 
\ref{dfn:B}. Then the image of the intersection of any orthant of 
$\Pro^{n+1}_\R$ with $\Xi_\cA$ under $\fii_\cA$ is either empty 
or a fixed affine hyperplane depending only on $\cA$. \qed 
\end{lemma} 

\noindent 
The proof of Lemma \ref{lemma:ckt} is simply a stream-lining of the proof of 
Theorem \ref{thm:gena}, based on the fact that the reduced $\cA$-discriminant 
contour is just a point when $\cA\!\in\!\R^{n\times (n+2)}$ has 
$\rank\!\left(\hA\right)\!=\!n+1$.  

\medskip 
\noindent 
{\bf Proof of Theorem \ref{thm:nearckt}:} The special case where 
$\conv\{a_1,\ldots,a_{n+3}\}$ is combinatorially simplicial 
follows immediately from Corollary \ref{cor:chamber} and 
Theorem \ref{thm:morse}. 

To address the general case, first assume  
$\cA$ is defective. So $\codim \Xi_\cA\!\geq\!2$, and thus 
(applying the Regular Interval Theorem again, as in Theorem 
\ref{thm:morse}) $\Pro^{t-1}_\R\setminus\Xi_\cA$ is path-connected, 
thereby implying that there is only $1$ isotopy type for $Z_\R(g)$. 
So, we may assume $\cA$ is non-defective. 

By Proposition \ref{prop:non}, we may also assume $\conv\{a_1,\ldots,a_{n+3}\}$ 
has exactly two facets containing exactly $n+1$ distinct points of 
$\{a_1,\ldots,a_{n+3}\}$. Proposition \ref{prop:2face} then tells us that the 
resulting signed reduce contour consists of $1$ locally convex 
arc $C$ (with $\leq\!n$ cusps) and $2$ lines we'll call 
$L_1$ and $L_2$. Considering the connected components of $\R^2\setminus C$, is 
easy to see that introducing $L_1$ creates no more than $n+2$ new 
connected components (since $L_1$ can intersect any 
sub-arc of $C$ in at most $2$ points, and $L_1$ is thus cut into 
at most $n+1$ sub-intervals). Similarly, counting the new intersections 
created by $(\R^2\setminus(C\cup L_1))\cap L_2$, we see that 
$\R^2\setminus(C\cup L_1\cup L_2)$ has no more than\\  
\mbox{}\hfill$\frac{1}{2}n^2-n+1+(n+2)+(n+3)\!=\!\frac{1}{2}n^2+\frac{3}{2}n+6$ 
\hfill\mbox{}\\  
connected components where $Z_\R(g)$ has constant topology 
for a fixed signed vector. So we are done. \qed 

\subsection{Proof of an Isotopy Lower Bound: Lemma \ref{lemma:weird}}  
\label{sub:weird} 
Substituting $x_i\!=\!e^{y_i}$, we then see that 
the lower bound for $\cA\!\in\!\Z^{1\times t}$ can be attained by counting the 
possible numbers of positive roots of a univariate $t$-nomial. Via a classic 
refinement of Descartes' Rule (see, e.g., \cite{descartes,grabiner}), a 
$t$-nomial can attain $k$ roots for any $k\!\in\!\{0,1,\ldots,t-1\}$ and we 
thus obtain the exact count of $t$ for the number of isotopy 
types of $g$. 

The lower bound for $\cA\!\in\!\Z^{2\times t}$ follows via Viro's Patchworking. 
In particular, in \cite{orevkov}, it is proved that the number of isotopy 
types for a 
real projective hypersurface defined by a homogeneous $(d+1)$-variate 
polynomial of degree $r$ is $2^{\Omega(r^d)}$. This construction in fact 
yields asymptotically the same number of isotopy types for a hypersurface 
in $\R^d_+$ defined by a $d$-variate polynomial of degree $r$, since 
we can use translation of the variables to move any ovals into the 
positive orthant. The number of monomial terms of such a polynomial 
is $\frac{(r+d)\cdots (r+1)}{d!}\!=\!\Theta(r^d)$ for fixed $d$ and increasing 
$r$. So this family of polynomials evinces $2^{\Omega(t)}$ distinct 
isotopy types for positive zero sets. Substituting exponentials for 
the variables, we are done. \qed 

\begin{rem} 
We have so far not addressed the conceptually simpler question of bounding the 
number of connected components of $Z_\R(g)$. In particular, we emphasize that 
the number of isotopy types can be far larger than the 
maximal number of connected components. For instance, in our preceding 
proof, we saw already that degree $r$ curves have $2^{\Omega(r^2)}$ 
isotopy types for their positive zero sets, while the maximal number of 
connected components of such a zero set is easily seen to be $O(r^2)$ 
by translating ovals and an application of Harnack's classical estimate.  
(\cite{viro16} contains an elegant discussion of Harnack's estimate.) \dia  
\end{rem} 

\section*{Acknowledgements}  
We thank Peter Gritzmann for important discussions on hyperplane 
arrangements. We also thank Mounir Nisse for insightful comments on 
a preliminary version of this work. This paper was developed during Rojas' 
recent visits to Technische Universit\"{a}t M\"{u}nchen and the Korean 
Institute of Advanced Studies. So we gratefully acknowledge the support of 
these institutions, and Rojas humbly thanks Peter Gritzmann and Mounir Nisse 
for their splendid hospitality. 

\bibliographystyle{amsalpha}

\end{document}